\newcommand{\xBc}{\langle}
\newcommand{\xBe}{\rangle}

\newcommand{\xbP}{\Pi}

\newcommand{\xbS}{\Sigma}

\newcommand{\xba}{\alpha}
\newcommand{\xbb}{\beta}

\newcommand{\xbe}{\in}
\newcommand{\xbf}{\phi}

\newcommand{\xbm}{\mu}

\newcommand{\xbo}{\omega}

\newcommand{\xbq}{\psi}
\newcommand{\xbr}{\rho}
\newcommand{\xbs}{\sigma}
\newcommand{\xbt}{\tau}

\newcommand{\xCK}{\times}

\newcommand{\xCN}{\neg}
\newcommand{\xCO}{ }
\newcommand{\xCQ}{\emptyset}

\newcommand{\xCf}{\hspace{0.1em}}

\newcommand{\xcA}{\forall}

\newcommand{\xcC}{\not\subseteq}

\newcommand{\xcE}{\exists}

\newcommand{\xcU}{\bigwedge}
\newcommand{\xcV}{\bigcup}

\newcommand{\xcb}{\subset}
\newcommand{\xcc}{\subseteq}

\newcommand{\xce}{\not\in}

\newcommand{\xcg}{\geq}
\newcommand{\xch}{\Rightarrow}
\newcommand{\xci}{\Leftarrow}

\newcommand{\xck}{\leq}

\newcommand{\xcm}{\models}

\newcommand{\xcp}{\rightarrow}

\newcommand{\xcr}{\leftrightarrow}
\newcommand{\xcs}{\cap}
\newcommand{\xcu}{\wedge}
\newcommand{\xcv}{\cup}

\newcommand{\xcz}{\Box}

\newcommand{\xDB}{\\[2ex]}

\newcommand{\xDH}{\item }

\newcommand{\xdR}{\Re}

\newcommand{\xda}{{\cal A}}

\newcommand{\xdp}{{\cal P}}

\newcommand{\xdx}{{\cal X}}
\newcommand{\xdy}{{\cal Y}}
\newcommand{\xdz}{{\cal Z}}

\newcommand{\xEH}{ & }
\newcommand{\xEI}{\begin{itemize}}
\newcommand{\xEJ}{\end{itemize}}
\newcommand{\xEP}{ \\ }

\newcommand{\xEd}{\neq}

\newcommand{\xEh}{\begin{enumerate}}

\newcommand{\xEj}{\end{enumerate}}

\newcommand{\xeB}{\not\prec}

\newcommand{\xeb}{\prec}
\newcommand{\xec}{\preceq}
\newcommand{\xed}{\succeq}

\newcommand{\xex}{\upharpoonright}

\newcommand{\xfA}{\mid}

\newcommand{\Xl}{\ldots}

\newcommand{\bl}{\begin{lemma} \rm}
\newcommand{\el}{\end{lemma}}
\newcommand{\br}{\begin{remark} \rm}
\newcommand{\er}{\end{remark}}
\newcommand{\be}{\begin{example} \rm}
\newcommand{\ee}{\end{example}}
\newcommand{\bco}{\begin{corollary} \rm}
\newcommand{\eco}{\end{corollary}}
\newcommand{\bc}{\begin{claim} \rm}
\newcommand{\ec}{\end{claim}}
\newcommand{\bfa}{\begin{fact} \rm}
\newcommand{\efa}{\end{fact}}
\newcommand{\bp}{\begin{proposition} \rm}
\newcommand{\ep}{\end{proposition}}
\newcommand{\bd}{\begin{definition} \rm}
\newcommand{\ed}{\end{definition}}
\newcommand{\bcs}{\begin{construction} \rm}
\newcommand{\ecs}{\end{construction}}
\newcommand{\bcd}{\begin{condition} \rm}
\newcommand{\ecd}{\end{condition}}
\newcommand{\bt}{\begin{theorem} \rm}
\newcommand{\et}{\end{theorem}}
\newcommand{\bn}{\begin{notation} \rm}
\newcommand{\en}{\end{notation}}
\newcommand{\bfi}{\begin{bild} \rm}
\newcommand{\efi}{\end{bild}}
\newcommand{\bsta}{\begin{statement} \rm}
\newcommand{\esta}{\end{statement}}
\newcommand{\bcom}{\begin{comment} \rm}
\newcommand{\ecom}{\end{comment}}
\newcommand{\bdia}{\begin{diagram} \rm}
\newcommand{\edia}{\end{diagram}}

\newcommand{\bfc}{\begin{figure}[htb] \begin{center}}
\newcommand{\efc}{\end{center} \end{figure}}

\documentclass{article}

\usepackage{amssymb,latexsym,epic,eepic,rotating}

\title{
Equilibria und weiteres Heiteres II-a
}

\author{Dov M. Gabbay
\thanks{
Dov.Gabbay@kcl.ac.uk, www.dcs.kcl.ac.uk/staff/dg
} \\
King's College, London
\thanks{
Department of Computer Science, King's College London, Strand,
London WC2R 2LS, UK
} \\
and \\
Bar-Ilan University, Israel
\thanks{
Department of Computer Science,
Bar-Ilan University,
52900 Ramat-Gan, Israel
} \\
and \\
University of Luxembourg
\thanks{
Computer Science and Communications,
Faculty of Sciences,
6, rue Coudenhove-Kalergi,
L-1359 Luxembourg
} \\ \\
Karl Schlechta
\thanks{
ks@cmi.univ-mrs.fr, karl.schlechta@web.de, http://www.cmi.univ-mrs.fr/ $\sim$ ks
} \\
Laboratoire d'Informatique Fondamentale de Marseille
\thanks{
CMI, 39, rue Joliot-Curie, F-13453 Marseille Cedex 13, France
(UMR 6166, CNRS and Universit\'{e} de Provence)
}
}

\begin{document}

\newtheorem{lemma}{Lemma}[section]
\newtheorem{theorem}[lemma]{Theorem}
\newtheorem{proposition}[lemma]{Proposition}
\newtheorem{corollary}[lemma]{Corollary}
\newtheorem{claim}[lemma]{Claim}
\newtheorem{fact}[lemma]{Fact}
\newtheorem{remark}[lemma]{Remark}
\newtheorem{definition}{Definition}[section]
\newtheorem{construction}{Construction}[section]
\newtheorem{condition}{Condition}[section]
\newtheorem{example}{Example}[section]
\newtheorem{notation}{Notation}[section]
\newtheorem{bild}{Figure}[section]
\newtheorem{comment}{Comment}[section]
\newtheorem{statement}{Statement}[section]
\newtheorem{diagram}{Diagram}[section]

\renewcommand{\labelenumi}
  {(\arabic{enumi})}
\renewcommand{\labelenumii}
  {(\arabic{enumi}.\arabic{enumii})}
\renewcommand{\labelenumiii}
  {(\arabic{enumi}.\arabic{enumii}.\arabic{enumiii})}
\renewcommand{\labelenumiv}
  {(\arabic{enumi}.\arabic{enumii}.\arabic{enumiii}.\arabic{enumiv})}

\maketitle
\setcounter{secnumdepth}{3}
\setcounter{tocdepth}{3}

\begin{abstract}

We investigate several technical and conceptual questions.

\end{abstract}

\tableofcontents

%
%
%
\section{
Introduction
}

We present here various results, which may one day be published in
a bigger paper, and which we wish to make already available to the
community.
\section{
Countably many disjoint sets
}

We show here that - independent of the cardinality of the language -
one can define only countably many inconsistent formulas.

The question is due to D. Makinson (personal communication).

$ \xCO $

We show here that, independent of the cardinality of the language,
one can define only countably many inconsistent formulas.

The problem is due to D. Makinson (personal communication).

\be

$\hspace{0.01em}$


\label{Example Co-Ex-Inf}

There is a countably infinite set of formulas s.t. the defined model sets
are pairwise disjoint.

Let $p_{i}:i \xbe \xbo $ be propositional variables.

Consider $ \xbf_{i}:= \xcU \{ \xCN p_{j}:j<i\} \xcu p_{i}$ for $i \xbe
\xbo.$

Obviously, $M(\xbf_{i}) \xEd \xCQ $ for all $i.$

Let $i<i';$ we show $M(\xbf_{i}) \xcs M(\xbf_{i' })= \xCQ.$ $M(
\xbf_{i' }) \xcm \xCN p_{i},$ $M(\xbf_{i}) \xcm p_{i}.$

$ \xcz $
\\[3ex]

\ee

\bfa

$\hspace{0.01em}$


\label{Fact Co-Ex-Inf}

Any set $X$ of consistent formulas with pairwise disjoint model sets is at
most
countable \index{countable}

\efa

\subparagraph{
Proof
}

$\hspace{0.01em}$


Let such $X$ be given.

(1) We may assume that $X$ consists of conjunctions of propositional
variables
or their negations.

Proof: Rewrite all $ \xbf \xbe X$ as disjunctions of conjunctions $
\xbf_{j}.$ At least one of
the conjunctions $ \xbf_{j}$ is consistent. Replace $ \xbf $ by one such $
\xbf_{j}.$ Consistency
is preserved, as is pairwise disjointness.

(2) Let $X$ be such a set of formulas. Let $X_{i} \xcc X$ be the set of
formulas in $X$ with
length $i,$ i.e., a consistent conjunction of $i$ many propositional
variables or
their negations, $i>0.$

As the model sets for $X$ are pairwise disjoint, the model sets for all $
\xbf \xbe X_{i}$
have to be disjoint.

(3) It suffices now to show that each $X_{i}$ is at most countable; we
even show
that each $X_{i}$ is finite.

Proof by induction:

Consider $i=1.$ Let $ \xbf, \xbf' \xbe X_{1}.$ Let $ \xbf $ be $p$ or $
\xCN p.$ If $ \xbf' $ is not $ \xCN \xbf,$ then
$ \xbf $ and $ \xbf' $ have a common model. So one must be $p,$ the other
$ \xCN p.$ But these
are all possibilities, so $card(X_{1})$ is finite.

Let the result be shown for $k<i.$

Consider now $X_{i}.$ Take arbitrary $ \xbf \xbe X_{i}.$ Without loss of
generality, let
$ \xbf =p_{1} \xcu  \Xl  \xcu p_{i}.$ Take arbitrary
$ \xbf' \xEd \xbf.$ As $M(\xbf) \xcs M(\xbf')= \xCQ,$ $ \xbf' $
must be a conjunction containing one of
$ \xCN p_{k},$ $1 \xck k \xck i.$ Consider now $X_{i,k}:=\{ \xbf' \xbe
X_{i}: \xbf' $ contains $ \xCN p_{k}\}.$
Thus $X_{i}=\{ \xbf \} \xcv \xcV \{X_{i,k}:1 \xck k \xck i\}.$ Note that
all $ \xbq, \xbq' \xbe X_{i,k}$ agree on $ \xCN p_{k},$
so the situation in $X_{i,k}$ is isomorphic to $X_{i-1}.$ So,
by induction hypothesis, $card(X_{i,k})$ is finite,
as all $ \xbf' \xbe X_{i,k}$ have to be mutually inconsistent. Thus,
$card(X_{i})$ is finite.
(Note that we did not use the fact that elements from different $X_{i,k},$
$X_{i,k' }$
also have to be mutually inconsistent; our rough proof suffices.)

$ \xcz $
\\[3ex]

Note that the proof depends very little on logic. We needed normal forms,
and used two truth values. Obviously, we can easily generalize to finitely
many truth values.

$ \xCO $
\section{
Independence as ternary relation
}

$ \xCO $
\subsection{
Introduction
}

\label{Section Prob-Func}
\subsubsection{
Independence
}

Independence is a central concept of reasoning.

In the context of non-monotonic logic and related areas like theory
revision,
it was perhaps first investigated formally by R. Parikh and co-authors, see
e.g.  \cite{Par96}, to
obtain ``local'' conflict solution.

The present authors investigated its role for interpolation in
preferential logics in
 \cite{GS10}, and showed connections to abstract multiplication of
size.

Independence plays also a central role for a FOL treatment of
preferential logics, where problems like the
``dark haired Swedes'' have to be treated. This is still subject
of ongoing research.

J. Pearl investigated independence in graphs and pobabilistic reasoning,
e.g. in  \cite{Pea88}, also as a ternary relation, $ \xBc X \xfA Y
\xfA Z \xBe.$

The aim of the present paper is to extend this abstract approach to the
preferential situation. We should emphasize that this is only an abstract
description of the independence relation, and thus not the same as
independence
for non-monotonic interpolation as examined in
 \cite{GS10}, where we $ \xCf used$ independence, essentially in the
form of the
multiplicative law
$ \xbm (X \xCK Y)= \xbm (X) \xCK \xbm (Y),$ which says that the $ \xbm
-$function preserves independence.

We have not investigated if an interesting form of interpolation results
from
some application of $ \xbm $ to situations described by $ \xBc X \xfA Y
\xfA Z \xBe,$ analogously
to above application of $ \xbm $ to situations described by $ \xBc X \xfA
\xfA Y \xBe.$
\subsubsection{
Overview
}

We will first discuss simple examples, to introduce the main ideas.

We then present the basic definitions formally, for probabilistic and
set independence.

We then show basic results for set independence as a ternary relation, and
turn
to our main
results, absence of finite characterization, and construction of new rules
for
this ternary relation.
\subsubsection{
Discussion of some simple examples
}

\label{Section Discussion-Simple}

We consider here $X=Y=Z=W=\{0,1\}$ and their products.
We will later generalize, but the main ideas stay the same.
First, we look at $X \xCK Z$ (the Cartesian product of $X$ with $Z),$ then
at
$X \xCK Z \xCK W,$ at $X \xCK Y \xCK Z,$ finally at $X \xCK Y \xCK Z \xCK
W.$ Elements of these products,
i.e., sequences, will be written for simplicity 00, 01, 10, etc., context
will disambiguate. General sequences will often be written $ \xbs,$ $
\xbt,$ etc.
We will also look at subsets of these products, like $\{00,11\} \xcc X
\xCK Z,$ and
various probability measures on these products.

As a matter of fact, the main part of this article concerns subsets A
of products $X_{1} \xCK  \Xl  \xCK X_{n}$ and a suitable notion of
independence for A,
roughly, if we can write A as $A_{1} \xCK  \Xl  \xCK A_{m}.$ This will be
made more
precise and discussed in progressively more complicated cases in this
section.

In the context of preferential structures, A is intended to be
$ \xbm (X_{1} \xCK  \Xl  \xCK X_{n}),$ the set of minimal models of $X_{1}
\xCK  \Xl  \xCK X_{n}.$
\paragraph{
$ X \xCK Z $
}

Let $P:X \xCK Z \xcp [0,1]$ be a (fixed) probability measure.

If $A \xcc X \xCK Z,$ we will set $P(A):= \xbS \{P(\xbs): \xbs \xbe
A\}.$

If $A_{x}:=\{ \xbs \xbe X \xCK Z: \xbs (X)=x\},$ we will write $P(x)$ for
$P(A_{x}),$ likewise
$P(z)$ for $P(A_{z}),$ if $A_{z}:=\{ \xbs \xbe X \xCK Z: \xbs (Z)=z\}.$
When these are ambiguous, we will e.g. write $A_{X=0}$ for $\{ \xbs \xbe X
\xCK Z: \xbs (X)=0\},$
and $P(X=0)$ for $P(A_{X=0}),$ etc.

We say that $X$ and $Z$ are independent for this $P$ iff for all
$xz \xbe X \xCK Z$ $P(xz)=P(x)*P(z).$

We write then $ \xBc X \xfA \xfA Z \xBe_{P},$
and call this and its variants probabilistic independence.

\be

$\hspace{0.01em}$


\label{Example XZ}

(1)

$P(00)=P(01)=1/6,$ $P(10)=P(11)=1/3.$

Then $P(X=0)=1/6+1/6=1/3,$ and $P(X=1)=2/3,$ $P(Z=0)=1/6+1/3=1/2,$ and
$P(Z=1)=1/2,$
so $ \xBc X \xfA \xfA Z \xBe_{P}.$

(2)

$P(00)=P(11)=1/3,$ $P(01)=P(10)=1/6.$

Then $P(X=0)=P(X=1)=P(Z=0)=P(Z=1)=1/2,$ but $P(00)=1/3 \xEd 1/2*1/2=1/4,$
so $ \xCN \xBc X \xfA \xfA Z \xBe_{P}.$

\ee

\bd

$\hspace{0.01em}$


\label{Definition P-A}

Consider now $ \xCQ \xEd A \xcc X \xCK Z$ for general $X,Z.$

Define the following probability measure on $X \xCK Z:$

\begin{flushleft}
\[  P_A(\xbs):=
\left\{ \begin{array}{lcl}
{\frac{1}{card(A)}}
\xEH iff \xEH  \xbs \xbe A  \xEP
\xEH \xEH \xEP
0 \xEH iff \xEH \xbs \xce A \xEP
\end{array}
\right.
\]
\end{flushleft}

\ed

\be

$\hspace{0.01em}$


\label{Example XZ-A}

(1)

$A:=\{00,01\},$

then $P_{A}(00)=P_{A}(01)=1/2,$ $P_{A}(10)=P_{A}(11)=0,$ $P_{A}(X=0)=1,$
$P_{A}(X=1)=0,$
$P_{A}(Z=0)=P_{A}(Z=1)=1/2,$ and we have $ \xBc X \xfA \xfA Z
\xBe_{P_{A}}.$

(2)

$A:=\{00,11\},$

then $P_{A}(00)=P_{A}(11)=1/2,$ $P_{A}(01)=P_{A}(10)=0,$
$P_{A}(X=0)=P_{A}(X=1)=1/2,$
$P_{A}(Z=0)=P_{A}(Z=1)=1/2,$ but $P_{A}(00)=1/2 \xEd
P_{A}(X=0)*P_{A}(Z=0)=1/4,$
and we have $ \xCN \xBc X \xfA \xfA Z \xBe_{P_{A}}.$

(3)

$A:=\{00,01,11\},$

then $P_{A}(00)=P_{A}(01)=P_{A}(11)=1/3,$ $P_{A}(10)=0,$ $P_{A}(X=0)=2/3,$
$P_{A}(X=1)=1/3,$
$P_{A}(Z=0)=1/3,$ $P_{A}(Z=1)=2/3,$
but $P_{A}(00)=1/3 \xEd P_{A}(X=0)*P_{A}(Z=0)=2/3*1/3=2/9,$
and we have $ \xCN \xBc X \xfA \xfA Z \xBe_{P_{A}}.$

\ee

Note that in (1) above, $A=\{0\} \xCK \{0,1\},$ but neither in (2), nor in
(3), A
can be written as such a product.
This is no coincidence, as we will see now.

More formally, we write $ \xBc X \xfA \xfA Z \xBe_{A}$ iff
for all $ \xbs \xbt \xbe A$ there is $ \xbr \xbe A$ such that $ \xbr (X)=
\xbs (X)$ and $ \xbr (Z)= \xbt (Z),$
or, equivalently, that $A=\{ \xbs (X): \xbs \xbe A\} \xCK \{ \xbs (Z):
\xbs \xbe A\},$
meaning that we can combine fragments of functions in A arbitrarily.

We call this and its variants set independence.

\bfa

$\hspace{0.01em}$


\label{Fact XZ-A}

Consider above situation $X \xCK Z.$ Then $ \xBc X \xfA \xfA Z
\xBe_{P_{A}}$ iff $ \xBc X \xfA \xfA Z \xBe_{A}.$

\efa

\subparagraph{
Proof
}

$\hspace{0.01em}$


``$ \xch $'':

$A \xcc \{ \xbs (X): \xbs \xbe A\} \xCK \{ \xbs (Z): \xbs \xbe A\}$ is
trivial.
Suppose $P_{A}(x,z)=P_{A}(x)*P_{A}(z),$ but there are $ \xbs, \xbt \xbe
A,$ $ \xbs (X) \xbt (Z) \xce A.$
Then $P_{A}(x),P_{A}(z)>0,$ but $P_{A}(x,z)=0,$ a contradiction.

``$ \xci $'':

Case 1: $P_{A}(x)=0,$ then $P_{A}(x,z)=0,$ and we are done. Likewise for
$P_{A}(Z)=0.$

Case 2: $P_{A}(x),P_{A}(z)>0.$

By definition and prerequisite,

$P_{A}(x)$ $=$ $ \frac{card\{ \xbs \xbe A: \xbs (X)=x\}}{card(A)}$ $=$ $
\frac{card\{ \xbs (Z): \xbs \xbe A\}}{card(A)},$

$P_{A}(z)$ $=$ $ \frac{card\{ \xbs \xbe A: \xbs (Z)=z\}}{card(A)}$ $=$ $
\frac{card\{ \xbs (X): \xbs \xbe A\}}{card(A)},$

$P_{A}(x,z)$ $=$ $ \frac{card\{ \xbs \xbe A: \xbs (X)=x, \xbs
(Z)=z\}}{card(A)}$ $=$ $ \frac{1}{card(A)}.$

By prerequisite again, $card(A)$ $=$ $card\{ \xbs (X): \xbs \xbe A\}$ $=$
$card\{ \xbs (Z): \xbs \xbe A\},$ so
$ \frac{card\{ \xbs (Z): \xbs \xbe A\}}{card(A)}$ $*$ $ \frac{card\{ \xbs
(X): \xbs \xbe A\}}{card(A)}$ $=$ $ \frac{1}{card(A)}$

$ \xcz $
\\[3ex]
\paragraph{
$ X \xCK Z \xCK W $
}

Here, $W$ will not be mentioned directly.

Let $P:X \xCK Z \xCK W \xcp [0,1]$ be a probability measure.

Again, we say that $X$ and $Z$ are independent for $P,$ $ \xBc X \xfA \xfA
Z \xBe_{P},$ iff
for all $x \xbe X,$ $z \xbe Z$ $P(x,z)=P(x)*P(z).$

\be

$\hspace{0.01em}$


\label{Example XZW}

(1)

Let $P(000)=P(001)=P(010)=P(011)=1/12,$ $P(100)=P(101)=P(110)=P(111)=1/6,$
then
$X$ and $Z$ are independent.

(2)

Let $P(100)=P(101)=P(010)=P(011)=1/12,$ $P(000)=P(001)=P(110)=P(111)=1/6,$
then
$P(X=0)=P(X=1)=P(Z=0)=P(Z=1)=1/2,$ but $P(X=0,Z=0)=1/3 \xEd 1/2*1/2=1/4,$
so $ \xCN \xBc X \xfA \xfA Z \xBe_{P}.$

\ee

As above, we define $P_{A}$ for $ \xCQ \xEd A \xcc X \xCK Z \xCK W.$

\be

$\hspace{0.01em}$


\label{Example XZW-A}

(1)

$A:=\{000,001,010,011\}.$ Then $P_{A}(X=0,Z=0)=P_{A}(X=0,Z=1)=1/2,$
$P_{A}(X=1,Z=0)=P_{A}(X=1,Z=1)=0,$ $P_{A}(X=0)=1,$ $P_{A}(X=1)=0,$
$P_{A}(Z=0)=P_{A}(Z=1)=1/2,$ so
$X$ and $Z$ are independent.

(2)

For $A:=\{000,001,110,111\},$ we see that $X$ and $Z$ are not independent
for $P_{A}.$

\ee

Considering possible decompositions of A into set products, we are not so
much
interested how many continuations into $W$ we have, but if there are any
or none.
This is often the case in logic, we are not interested how many models
there
are, but if there is a model at all.

Thus we define independence for A again by:

$ \xBc X \xfA \xfA Z \xBe_{A}$ iff
for all $ \xbs \xbt \xbe A$ there is $ \xbr \xbe A$ such that $ \xbr (X)=
\xbs (X)$ and $ \xbr (Z)= \xbt (Z).$

The equivalence between probabilitistic independence,
$ \xBc X \xfA \xfA Z \xBe_{P_{A}}$ and set independence, $ \xBc X \xfA
\xfA Z \xBe_{A}$ is lost now, as the
second part of the
following example shows:

\be

$\hspace{0.01em}$


\label{Example XZW-A-2}

(1)

$A:=\{000,010,100,110\}$ satisfies both forms of independence,
$ \xBc X \xfA \xfA Z \xBe_{P_{A}}$ and set independence, $ \xBc X \xfA
\xfA Z \xBe_{A}.$

\ee

(2)

$A:=\{000,001,010,100,110\}.$

Here, we have $P_{A}(X=0)=3/5,$ $P_{A}(X=1)=2/5,$ $P_{A}(Z=0)=3/5,$
$P_{A}(Z=1)=2/5,$
but $P_{A}(X=0,Z=0)=2/5 \xEd 3/5*3/5.$

Consider now $ \xBc X \xfA \xfA Z \xBe_{A}:$ Take $ \xbs, \xbt \xbe A,$
then for all possible values
$ \xbs (X),$ $ \xbt (Z),$ there is $ \xbr $ such that $ \xbr (X)= \xbs
(X),$ $ \xbr (Z)= \xbt (Z)$ - the value
$ \xbr (W)$ is without importance.

We have, however:

\bfa

$\hspace{0.01em}$


\label{Fact XZW-A}

$ \xBc X \xfA \xfA Z \xBe_{P_{A}}$ $ \xch $ $ \xBc X \xfA \xfA Z
\xBe_{A}.$

\efa

\subparagraph{
Proof
}

$\hspace{0.01em}$


Let $ \xbs, \xbt \xbe A,$ but suppose there is no $ \xbr \xbe A$ such
that $ \xbr (X)= \xbs (X)$ and
$ \xbr (Z)= \xbt (Z).$ Then $P_{A}(\xbs (X)),P_{A}(\xbt (Z))>0,$ but
$P_{A}(\xbs (X), \xbt (Z))=0.$ $ \xcz $
\\[3ex]
\paragraph{
$ X \xCK Y \xCK Z $
}

\label{Section XYZ}

We consider now independence of $X$ and $Z,$ given $Y.$

The probabilistic definition is:

$ \xBc X \xfA Y \xfA Z \xBe_{P}$ iff for all $x \xbe X,y \xbe Y,z \xbe Z$
$P(x,y,z)*P(y)=P(x,y)*P(y,z).$

As we are interested mainly in subsets $A \xcc X \xCK Y \xCK Z$ and the
resulting $P_{A},$
and combination of function fragments, we work immediately with these.

We have to define $ \xBc X \xfA Y \xfA Z \xBe_{A}.$

$ \xBc X \xfA Y \xfA Z \xBe_{A}$ iff for all $ \xbs, \xbt \xbe A$ such
that $ \xbs (Y)= \xbt (Y)$ there is $ \xbr \xbe A$
such that $ \xbr (X)= \xbs (X),$ $ \xbr (Y)= \xbs (Y)= \xbt (Y),$ $ \xbr
(Z)= \xbt (Z).$

When we set for $y \xbe Y$ $A_{y}:=\{ \xbs \xbe A: \xbs (Y)=y\},$ we then
have:

$A_{y}=\{ \xbs (X): \xbs \xbe A_{y}\} \xCK \{y\} \xCK \{ \xbs (Z): \xbs
\xbe A_{y}\}.$

The following example shows that $ \xBc X \xfA Y \xfA Z \xBe_{A}$ and $
\xBc X \xfA \xfA Z \xBe_{A}$ are independent
from each other:

\be

$\hspace{0.01em}$


\label{Example XYZ-A}

(1)

$ \xBc X \xfA Y \xfA Z \xBe_{A}$ may hold, but not $ \xBc X \xfA \xfA Z
\xBe_{A}:$

Consider $A:=\{000,111\}.$ $ \xBc X \xfA Y \xfA Z \xBe_{A}$ is obvious, as
only $ \xbs $ goes through
each element in the middle. But there is no 0x1, so $ \xBc X \xfA \xfA Z
\xBe_{A}$ fails.

(2)

$ \xBc X \xfA \xfA Z \xBe_{A}$ may hold, but not $ \xBc X \xfA Y \xfA Z
\xBe_{A}:$

Consider $A:=\{000,101,110,011\}.$ Fixing, e.g., 0 in the middle shows
that
$ \xBc X \xfA Y \xfA Z \xBe_{A}$ fails, but neglecting the middle, we can
combine arbitrarily, so
$ \xBc X \xfA \xfA Z \xBe_{A}$ holds.

\ee

\be

$\hspace{0.01em}$


\label{Example XYZ-Prod}

This example show that $ \xBc X \xfA Y \xfA Z \xBe_{A}$ does not mean that
A is some
product $A_{X} \xCK A_{Y} \xCK A_{Z}:$

Let $A:=\{000,111\},$ then clearly $ \xBc X \xfA Y \xfA Z \xBe_{A},$ but A
is no such product.

\ee

We have again:

\bfa

$\hspace{0.01em}$


\label{Fact XYZ-A}

Let $ \xCQ \xEd A \xcc X \xCK Y \xCK Z,$ then
$ \xBc X \xfA Y \xfA Z \xBe_{A}$ and $ \xBc X \xfA Y \xfA Z \xBe_{P_{A}}$
are equivalent.

\efa

\subparagraph{
Proof
}

$\hspace{0.01em}$


``$ \xci $'':

Suppose there are $ \xbs, \xbt \xbe A$ such that $ \xbs (Y)= \xbt (Y),$
but there is no $ \xbr \xbe A$
such that $ \xbr (X)= \xbs (X),$ $ \xbr (Y)= \xbs (Y)= \xbt (Y),$ $ \xbr
(Z)= \xbt (Z).$
Then $P_{A}(\xbs (X), \xbs (Y)),P_{A}(\xbt (Y), \xbt (Z)),P_{A}(\xbs
(Y))>0,$ but
$P_{A}(\xbs (X), \xbs (Y)= \xbt (Y), \xbt (Z))=0.$

``$ \xch $'':

Case 1: $P_{A}(x,y)$ or $P_{A}(y,z)=0,$ then $P_{A}(x,y,z)=0,$ and we are
done.

Case 2: $P_{A}(x,y),P_{A}(y,z)>0.$
By definition and prerequisite,
$P_{A}(x,y)$ $=$ $ \frac{card\{ \xbs \xbe A: \xbs (X)=x, \xbs
(Y)=y\}}{card(A)}$ $=$ $ \frac{card\{ \xbs (Z): \xbs \xbe A, \xbs
(Y)=y\}}{card(A)}$ and
$P_{A}(y,z)$ $=$ $ \frac{card\{ \xbs \xbe A: \xbs (Y)=y, \xbs
(Z)=Z\}}{card(A)}$ $=$ $ \frac{card\{ \xbs (X): \xbs \xbe A, \xbs
(Y)=y\}}{card(A)},$ so
$P_{A}(x,y)*P_{A}(y,z)$ $=$ $ \frac{card\{ \xbs \xbe A: \xbs
(Y)=y\}}{card(A)*card(A)}.$ Moreover,
$P_{A}(y)$ $=$ $ \frac{card\{ \xbs \xbe A: \xbs (Y)=y\}}{card(A)},$
$P_{A}(x,y,z)$ $=$ $ \frac{1}{card(A)},$ so
$P_{A}(y)*P_{A}(x,y,z)$ $=$ $ \frac{card\{ \xbs \xbe A: \xbs
(Y)=y\}}{card(A)*card(A)}$ $=$ $P_{A}(x,y)*P_{A}(y,z)$

$ \xcz $
\\[3ex]
\paragraph{
$ X \xCK Y \xCK Z \xCK W $
}

The definitions stay the same as for $X \xCK Y \xCK Z.$

The equivalence between probabilitistic independence,
$ \xBc X \xfA Y \xfA Z \xBe_{P_{A}}$ and set independence, $ \xBc X \xfA Y
\xfA Z \xBe_{A}$ is lost again, as the
following example shows:

\be

$\hspace{0.01em}$


\label{Example XYZW-A}

$A:=\{0000,0001,0010,1000,1010\}.$

Here, we have $P_{A}(X=0,Y=0)=3/5,$ $P_{A}(X=1,Y=0)=2/5,$
$P_{A}(Y=0,Z=0)=3/5,$ $P_{A}(Y=0,Z=1)=2/5,$ $P_{A}(Y=0)=1,$
but $P_{A}(X=0,Y=0,Z=0)=2/5 \xEd 3/5*3/5.$

Consider now $ \xBc X \xfA Y \xfA Z \xBe_{A}:$ Take $ \xbs, \xbt \xbe A,$
such that $ \xbs (Y)= \xbt (Y),$
then for all possible values
$ \xbs (X),$ $ \xbt (Z),$ there is $ \xbr $ such that $ \xbr (X)= \xbs
(X),$ $ \xbr (Y)= \xbs (Y)= \xbt (Y),$
$ \xbr (Z)= \xbt (Z)$ - the value $ \xbr (W)$ is without importance.

\ee

We have, however:

\bfa

$\hspace{0.01em}$


\label{Fact XYZW-A}

$ \xBc X \xfA Y \xfA Z \xBe_{P_{A}}$ $ \xch $ $ \xBc X \xfA Y \xfA Z
\xBe_{A}.$

\efa

\subparagraph{
Proof
}

$\hspace{0.01em}$


Let $ \xbs, \xbt \xbe A$ such that $ \xbs (Y)= \xbt (Y),$ but suppose
there is no $ \xbr \xbe A$ such that
$ \xbr (X)= \xbs (X),$ $ \xbr (Y)= \xbs (Y)= \xbt (Y),$
$ \xbr (Z)= \xbt (Z).$ Then $P_{A}(\xbs (X), \xbs (Y)),P_{A}(\xbs (Y),
\xbt (Z))>0,$
but $P_{A}(\xbs (X), \xbs (Y), \xbt (Z))=0.$ $ \xcz $
\\[3ex]
\paragraph{
A remark on generalization
}

The $X,Y,Z,W$ may also be more complicated sets, themselves products,
but this will not change definitions and results beyond notation.

In the more complicated cases, we will often denote subsets by more
complicated letters than A, e.g., by $ \xbS.$
\paragraph{
A remark on intuition
}

Consider set independence, where $A:= \xbm (U),$ $U=U_{1} \xCK  \Xl  \xCK
U_{n}.$
Set $ \xBc  \Xl  \xBe:= \xBc  \Xl  \xBe_{ \xbm (U)}.$

 \xEh

 \xDH
$ \xBc X \xfA \xfA Z \xBe $ means then:

 \xEh
 \xDH
all we know is that we are in a normal situation,
 \xDH
if we know in addition something definite about $Z$ (1 model!)
we do not know anything more about $X,$ and vice versa.
 \xEj

$ \xBc X \xfA Y \xfA Z \xBe $ means then:

 \xEh
 \xDH
all we know is that we are in a normal situation,
 \xDH
if we have definite information about $Y,$
we may know more about $X.$ But knowing something
in addition about $Z$ will not give us not more information about $X,$ and
conversely.
 \xEj

 \xDH
The restriction to $ \xbm (U)$ codes our background knowledge.

 \xDH
Note that $X \xcv Y \xcv Z$ need not be $I,$ e.g., $W$ might be missing.
We did not
count the continuations into $W,$ but considered only existence of a
continuation
(if this does not exist, then there just is no such sequence).

This corrsponds to multiplication with 1, the unit ALL on $W,$ or, more
generally, in the rest of the paper, with $1_{I-(X \xcv Y \xcv Z)}.$
We may choose however
we want, it has to be somewhere, in ALL.

 \xEj
\subsubsection{
Basic definitions
}

\bd

$\hspace{0.01em}$


\label{Definition Restrict}

If $f$ is a function, $Y$ a subset of its domain, we write
$f \xex Y$ for the restriction of $f$ to elements of $Y.$

If $F$ is a set of functions over $Y,$ then $F \xex Y:=\{f \xex Y:f \xbe
F\}.$
\subsection{
Probabilistic and set independence
}
\subsubsection{
Probabilistic independence
}

\ed

Independence as an abstract ternary relation for probability and other
situations has been examined by W. Spohn,
see  \cite{Spo80}, A. P. Dawid,
see  \cite{Daw79},
J. Pearl, see, e.g.,
 \cite{Pea88}, etc.

\bd

$\hspace{0.01em}$


\label{Definition Restrict-2}

(1)

Let $I \xEd \xCQ $ be an arbitrary (index) set, for $i \xbe I$ $U_{i} \xEd
\xCQ $ arbitrary sets.
Let $U:= \xbP \{U_{i}:i \xbe I\},$ and for $X \xcc I$ $U_{X}:= \xbP
\{U_{i}:i \xbe X\}.$

(2)

Let $P: \xdp (U) \xcp [0,1]$ be a probability measure. (We may assume that
$P$ is defined
by its value on singletons.)

(3.1)

By abuse of language, for $X \xcc I,$ $x \xbe U_{X},$
let $P(x)$ $:=$ $P(\{u \xbe U: \xcA i \xbe Xu(i)=x(i)\}),$ so
$P(x)=P(\{u \xbe U:u \xex X=x\}).$

Analogously, for $X,Y \xcc I,$ $X \xcs Y= \xCQ,$ $x \xbe U_{X},$ $y \xbe
U_{Y},$
let $P(x,y)$ $:=$ $P(\{u \xbe U:$ $u \xex X=x$ and $u \xex Y=y\}).$

(3.2)

Finally, for $X,Y,Z \xcc I$ pairwise disjoint, $x \xbe U_{X},$ $y \xbe
U_{Y},$ $z \xbe U_{Z},$ let
$P(x \xfA y):= \frac{P(x,y)}{P(y)},$ $P(x \xfA y,z):=
\frac{P(x,y,z)}{P(y,z)},$ etc.

(We have, of course, to pay attention that we do not divide by 0.)

\ed

\bd

$\hspace{0.01em}$


\label{Definition Prob-Ind}

$P$ as above defines a 3-place relation of independence
on pairwise disjoint $X,Y,Z \xcc I$
$ \xBc X \xfA Y \xfA Z \xBe_{P}$ by

\begin{flushleft}
\[  \xBc  \xdx \xfA \xdy \xfA \xdz  \xBe _{P}:\xcr
\left\{ \begin{array}{lcl}
 \xcA x \xbe U_X,  \xcA y \xbe U_Y,  \xcA z \xbe U_Z
(P(y,z)> 0  \xcp  P(x \xfA y) =P(x \xfA y,z)),
\xEH if \xEH  Y \xEd \xCQ  \xEP
i.e., P(x,y)/P(y)=P(x,y,z)/P(y,z), or
\xEH \xEH \xEP
P(x,y,z)*P(y)=P(x,y)*P(y,z)
\xEH \xEH \xEP
\xEH \xEH \xEP
\xEH \xEH \xEP
 \xcA x \xbe U_X,  \xcA z \xbe U_Z  (P(z)>0  \xcp
P(x)= P(x \xfA z)),
\xEH if \xEH  Y  = \xCQ  \xEP
i.e., P(x)=P(x,z)/P(z), or
\xEH \xEH \xEP
P(x,z)=P(x)*P(z)
\xEH \xEH \xEP
\end{array}
\right.
\]
\end{flushleft}

If $Y= \xCQ,$ we shall also write $ \xBc X \xfA \xfA Z \xBe_{P}$ for $
\xBc X \xfA Y \xfA Z \xBe_{P}$.

Recall from
Section 
\ref{Section Discussion-Simple} (page 
\pageref{Section Discussion-Simple})  that we call
this notion probabilistic independence.

\ed

E.g., Pearl discusses the rules $(a)-(e)$ of
Definition 
\ref{Definition Basic-Rules} (page 
\pageref{Definition Basic-Rules})  for the relation
defined in Definition 
\ref{Definition Prob-Ind} (page 
\pageref{Definition Prob-Ind}).

\bd

$\hspace{0.01em}$


\label{Definition Basic-Rules}

(a) Symmetry: $ \xBc X \xfA Y \xfA Z \xBe $ $ \xcr $ $ \xBc Z \xfA Y \xfA
X \xBe $

(b) Decomposition: $ \xBc X \xfA Y \xfA Z \xcv W \xBe $ $ \xcp $ $ \xBc X
\xfA Y \xfA Z \xBe $

(c) Weak Union: $ \xBc X \xfA Y \xfA Z \xcv W \xBe $ $ \xcp $ $ \xBc X
\xfA Y \xcv W \xfA Z \xBe $

(d) Contraction: $ \xBc X \xfA Y \xfA Z \xBe $ and $ \xBc X \xfA Y \xcv Z
\xfA W \xBe $ $ \xcp $ $ \xBc X \xfA Y \xfA Z \xcv W \xBe $

(e) Intersection: $ \xBc X \xfA Y \xcv W \xfA Z \xBe $ and $ \xBc X \xfA
Y \xcv Z \xfA W \xBe $ $ \xcp $ $ \xBc X \xfA Y \xfA Z \xcv W \xBe $

$(\xCQ)$ Empty outside: $ \xBc X \xfA Y \xfA Z \xBe $ if $X= \xCQ $ or
$Z= \xCQ.$

\ed

\bp

$\hspace{0.01em}$


\label{Proposition Prob-Valid}

If $P$ is a probability measure, and $ \xBc X \xfA Y \xfA Z \xBe_{P}$
defined as above, then
$(a)-(d)$ of Definition 
\ref{Definition Basic-Rules} (page 
\pageref{Definition Basic-Rules})  hold for $ \xBc  \Xl
\xBe = \xBc  \Xl  \xBe_{P},$
and if $P$ is strictly positive,
(e) will also hold.

\ep

The proof is elementary, well known, and will not be repeated here.

Doch ein Beispiel geben?
\paragraph{
A side remark on preferential structures
}

Being a minimal element is not upward absolute in general preferential
structures, but in raked structures, provided the smaller set contains
some element minimal in the bigger set.

\bfa

$\hspace{0.01em}$


\label{Fact Basic-Pro}

In the probabilistic interpretation, the following holds:

Let $U$ be a finite set, $f:U \xcp \xdR $ such that $ \xcA u \xbe U.f(u)
\xcg 0.$

For all $A \xcc U,$ such that $ \xcE a' \xbe A.f(a')>0$ and all $a \xbe
A$

$f_{A}(a):= \frac{f(a)}{ \xbS \{f(a'):a' \xbe A\}}$ defines a probability
measure on $ \xCf A.$

For $B \xcc A,$ define $f_{A}(B):= \xbS \{f_{A}(b):b \xbe B\}.$ Then the
following property
holds:

(BASIC) For all $D \xcc B \xcc A \xcc U$ such that $ \xcE b \xbe B.f(b)>0$
$f_{A}(D)=f_{A}(B)*f_{B}(D).$

\efa

\subparagraph{
Proof
}

$\hspace{0.01em}$


For $X \xcc Y \xcc U$ such that $ \xcE y \xbe Y.f(y)>0$ we have
$f_{Y}(X):= \xbS \{f_{Y}(x):x \xbe X\}= \frac{ \xbS \{f(x):x \xbe X\}}{
\xbS \{f(y):y \xbe Y\}}.$

Thus, $f_{A}(D)$ $:=$ $ \frac{ \xbS \{f(d):d \xbe D\}}{ \xbS \{f(a):a \xbe
A\}}$ $=$ $ \frac{ \xbS \{f(b):b \xbe B\}}{ \xbS \{f(a):a \xbe A\}}$ $*$ $
\frac{ \xbS \{f(d):d \xbe D\}}{ \xbS \{f(b):b \xbe
B\}}$ $=$ $f_{A}(B)*f_{B}(D).$

$ \xcz $
\\[3ex]

We have the following fact for $ \xbm $ generated by a relation:

\bfa

$\hspace{0.01em}$


\label{Fact Basic-Pref}

Let $U$ be a finite preferential structure such that for $A \xcc U$ $ \xbm
(A)= \xCQ $ $ \xch $ $A= \xCQ.$

Then $U$ is ranked iff (BASIC) as defined in
Fact \ref{Fact Basic-Pro} (page \pageref{Fact Basic-Pro})  holds for $f_{A}.$

\efa

\subparagraph{
Proof
}

$\hspace{0.01em}$


``$ \xch $'':

Let $D \xcc B \xcc A \xcc U,$ $B \xEd \xCQ.$

Case 1: $D \xcs \xbm (A)= \xCQ.$ Then $f_{A}(D)=0.$

Case 1.1: If $B \xcs \xbm (A)= \xCQ,$ then $f_{A}(B)=0,$ and we are done.

Case 1.2: Let $B \xcs \xbm (A) \xEd \xCQ.$ If $D \xcs \xbm (B)= \xCQ,$
then $f_{B}(D)=0,$ and we are done.
Suppose $D \xcs \xbm (B) \xEd \xCQ,$ so there is $d \xbe D \xcs \xbm
(B),$
so $d \xbe D \xcs \xbm (A)$ by $B \xcs \xbm (A) \xEd \xCQ $ and
rankedness, so $f_{A}(D) \xEd \xCQ,$
contradiction.

Case 2: $D \xcs \xbm (A) \xEd \xCQ.$

Thus, by $D \xcc B,$ $B \xcs \xbm (A) \xEd \xCQ,$ and by rankedness $
\xbm (B)=B \xcs \xbm (A).$ So
by $D \xcc B$ again, $D \xcs \xbm (A)=D \xcs (B \xcs \xbm (A))=D \xcs \xbm
(B).$
By definition,
$f_{A}(B):= \frac{card(\xbm (A) \xcs B)}{card(\xbm (A))},$
$f_{A}(D):= \frac{card(\xbm (A) \xcs D)}{card(\xbm (A))},$
$f_{B}(D):= \frac{card(\xbm (B) \xcs D)}{card(\xbm (B))}.$
Thus,
$ \frac{card(\xbm (A) \xcs D)}{card(\xbm (A))}= \frac{card(\xbm (A)
\xcs B)}{card(\xbm (A))}* \frac{card(\xbm (B) \xcs D)}{card(\xbm
(B))}.$

``$ \xci $'':

Then there are $a,b,c \xbe U,$ where $ \xCf a$ is incomparable to $b,$ and
$b \xeb c$ but $a \xeB c,$
or $c \xeb b,$ but $c \xeB a.$ We have four possible cases.

Let, in all cases, $A:=\{a,b,c\}.$ We construct a contradiction to
(BASIC).

Case 1, $b \xeb c:$

Case 1.1, a is incomparable to $c:$
Consider $B:=\{a,c\},$ $D:=\{a\}.$
Then $f_{A}(D)= \frac{1}{2},$ $f_{A}(B)= \frac{1}{2},$ $f_{B}(D)=
\frac{1}{2}.$

Case 1.2, $c \xeb a$ (so $ \xeb $ is not transitive):
Consider $B:=\{a,b\},$ $D:=\{a\}.$
Then $f_{A}(D)=0,$ $f_{A}(B)=1,$ $f_{B}(D)= \frac{1}{2}.$

Case 2, $c \xeb b:$

Case 2.1, a is incomparable to $c:$

Consider $B:=\{a,b\},$ $D:=\{a\}.$ Then
$f_{A}(D)= \frac{1}{2},$ $f_{A}(B)= \frac{1}{2},$ $f_{B}(D)= \frac{1}{2}.$

Case 2.2, $a \xeb c$ - similar to Case 1.2.

$ \xcz $
\\[3ex]

\br

$\hspace{0.01em}$


\label{Remark Zero}

Note that sets $A \xcc B,$ where $ \xbm (B) \xcs A= \xCQ,$ and sets where
$P(A)=0$ have a similar,
exceptional role. This might still be important.
\subsubsection{
Set independence
}

\er

We interpret independence here differently, but in a related way,
as prepared in
Section 
\ref{Section Discussion-Simple} (page 
\pageref{Section Discussion-Simple}).

\bd

$\hspace{0.01em}$


\label{Definition Indep-Fct}

We consider function sets $ \xbS $ etc. over a fixed, arbitrary domain $I
\xEd \xCQ,$ into
some fixed codomain $K.$

(1)

For pairwise disjoint subsets $X,Y,Z$ of $I,$ we define

$ \xBc X \xfA Y \xfA Z \xBe _{ \xbS }$ iff for all $f,g \xbe \xbS $ such that $f
\xex Y=g \xex Y,$
there is $h \xbe \xbS $ such that $h \xex X=f \xex X,$ $h \xex Y=f \xex
Y=g \xex Y,$
$h \xex Z=g \xex Z.$

Recall from
Section 
\ref{Section Discussion-Simple} (page 
\pageref{Section Discussion-Simple})  that we call
this notion set independence.

$Y$ may be empty, then the condition $f \xex Y=g \xex Y$ is void.

Note that nothing is said about $I-(X \xcv Y \xcv Z),$ so we look at the
projection
of $U$ to $X \xcv Y \xcv Z.$

When $Y= \xCQ,$ we will also write $ \xBc X \xfA \xfA Z \xBe _{ \xbS }.$

$ \xBc X \xfA Y \xfA Z \xBe _{ \xbS }$ means thus, that we can piece functions
together, or that we have a
sort of decomposition of $ \xbS $ into a product. This is an independence
property,
we can put parts together independently.

(2)

In the sequel, we will just write $ \xBc  \Xl  \xBe $ for $ \xBc  \Xl
\xBe_{ \xbS }$ when the meaning is
clear from the context.

\ed

Recall that
Example 
\ref{Example XZW-A-2} (page 
\pageref{Example XZW-A-2})  compares different forms of independence,
the
probabilistic and the set variant.

Obviously, we can generalize the equivalence results for probabilistic
and set independence for $X \xCK Z$ and $X \xCK Y \xCK Z$ to the general
situation with $W$ in
Section 
\ref{Section Discussion-Simple} (page 
\pageref{Section Discussion-Simple}),
as long as we do not consider the full functions $ \xbs,$ but only their
restrictions to $X,Y,Z,$ $ \xbs \xex (X \xcv Y \xcv Z).$
As we will stop the discussion of probablistic independece here, and
restrict ourselves to set independence, this is left as an easy exercise
to the reader.
\clearpage
\subsection{
Basic results for set independence
}

\bn

$\hspace{0.01em}$


\label{Notation Indep-Fct}

In more complicated cases,
we will often write $ \xCf ABC$ for $ \xBc A \xfA B \xfA C \xBe,$ and $ \xCN
ABC$
or $- \xCf ABC$ if $ \xBc A \xfA B \xfA C \xBe $ does
not hold.
Moreover, we will often just write $f(A)$ for $f \xex A,$ etc.

For $ \xBc A \xcv A' \xfA B \xfA C \xBe,$ we will then write $(AA')BC,$ etc.

If only singletons are involved, we will sometimes write $ \xCf abc$
instead of $ \xCf ABC,$
etc.

When we speak about fragments of functions, we will often write just
$A: \xbs $ for $ \xbs \xex A,$ $B: \xbs = \xbt $ for $ \xbs \xex B= \xbt
\xex B,$ etc.

\en

We use the following notations for functions:

\bd

$\hspace{0.01em}$


\label{Definition Const-Func-H}

The constant functions $0_{c}$ and $1_{c}:$

$0_{c}(i)=0$ for all $i \xbe I$

$1_{c}(i)=1$ for all $i \xbe I$

Moreover, when we define a function $ \xbs:I \xcp \{0,1\}$ argument by
argument,
we abbreviate $ \xbs (a)=0$ by $a=0,$ etc.

Sometimes, we also give (a fragment of) a function just by the sequence of
the
values, so instead of writing $a=0,$ $b=1,$ $c=1,$ we just write 011 -
context
will disambiguate.

\ed

\br

$\hspace{0.01em}$


\label{Remark Sys-Valid-H}

This remark gives an intuitive justification of (some of) above rules in
our context.

Rule (a) is trivial.

It is easiest to set $Y:= \xCQ $ to see the intuitive meaning.

Rule (b) is a trivial consequence. If we can combine longer sequences,
then we
can combine shorter, too.

Rule (c) is again a trivial consequence. If we can combine arbitrary
sequences,
then we can also combine those which agree already on some part.

Rule (d) is the most interesting one, it says when we may combine $ \xCf
longer$
sequences. Having just $ \xBc X \xfA \xfA Z \xBe $ and $ \xBc X \xfA \xfA
W \xBe $ as prerequisite does not
suffice, as we might lose when applying $ \xBc X \xfA \xfA W \xBe $ what
we had already by
$ \xBc X \xfA \xfA Z \xBe.$ The condition $ \xBc X \xfA Z \xfA W \xBe $
guarantees that we do not lose this.

In our context, it means the following:

We want to combine $ \xbs \xex X$ with $ \xbt \xex Z \xcv W.$
By $ \xBc X \xfA \xfA Z \xBe,$ we can combine $ \xbs \xex X$ with $ \xbt
\xex Z.$ Fix $ \xbr $ such that
$ \xbr \xex X= \xbs \xex X,$ $ \xbr \xex Z= \xbt \xex Z.$ As $ \xbr \xex
Z= \xbt \xex Z,$ by $ \xBc X \xfA Z \xfA W \xBe,$ we can
combine $ \xbr \xex X \xcv Z$ with $ \xbt \xex W,$ and have the result.

Note that we change the functions here, too: we start with $ \xbs,$ $
\xbt,$ then continue
with $ \xbr,$ $ \xbt.$

We can use what we constructed already as a sort of scaffolding for
constructing the rest.

\er

\bfa

$\hspace{0.01em}$


\label{Fact Sys-Prod}

Zusammenhang $ \xBc X \xfA Y \xfA Z \xBe $ mit Produkten.

\efa

\subparagraph{
Proof
}

$\hspace{0.01em}$


Do

$ \xcz $
\\[3ex]

We show now that above Rules $(a)-(d)$ hold in our context, but (e) does
not hold.

\bfa

$\hspace{0.01em}$


\label{Fact Sys-Valid-H}

In our interpretation,

(1)
rule (e) does not hold,

(2)
all $ \xBc X \xfA Y \xfA \xCQ \xBe $ (and thus also all $ \xBc \xCQ \xfA Y
\xfA Z \xBe)$ hold.

(3)
rules $(a)-(d)$ hold, even when one or both of the outside elements of the
tripels is the empty set.

\efa

\subparagraph{
Proof
}

$\hspace{0.01em}$


(1) (e) does not hold:

Consider $I:=\{x,y,z,w\}$ and
$U:=\{1111,0100\}.$ Then $x(yw)z$ and $x(yz)w,$
as for all $ \xbs \xex yw$ there is just one $ \xbt $
this $ \xbs $ can be. The same holds for $x(yz)w.$
But for $y=1,$ there are two different paths through $y=1,$ which cannot
be
combined.

(2) This is a trivial consequence of the fact that $\{f:$ $f: \xCQ \xcp
U\}=\{ \xCQ \}.$

(3)
Rules (a), (b), (c) are trivial, by definition, also for $X,Z= \xCQ.$
In (c), if $W= \xCQ,$ there is nothing to show.

Rule (d): The cases for $X,W,Z= \xCQ $ are trivial.
Assume $ \xbs,$ $ \xbt $ such that $ \xbs \xex Y= \xbt \xex Y,$ we want
to combine $ \xbs \xex X$ with
$ \xbt \xex Z \xcv W.$ By $ \xBc X \xfA Y \xfA Z \xBe,$ there is $ \xbr $
such that $ \xbr \xex X= \xbs \xex X,$ $ \xbr \xex Y= \xbs \xex Y= \xbt
\xex Y,$
$ \xba \xex X= \xbr \xex Z= \xbt \xex Z.$ Thus $ \xbr $ and $ \xbt $
satisfy the prerequisite of
$ \xBc X \xfA Y \xcv Z \xfA W \xBe,$ and there is $ \xba $ such that $
\xba \xex X= \xbr \xex X= \xbs \xex X,$
$ \xba \xex X= \xbr \xex Y= \xbs \xex Y= \xbt \xex Y,$ $ \xba \xex W= \xbt
\xex W.$

$ \xcz $
\\[3ex]

Next, we give examples which shows that increasing the center set
can change validity of the tripel in any way.

\be

$\hspace{0.01em}$


\label{Example Change-Rel}

(1)

This example shows that neither $ \xBc X \xfA Y \xfA Z \xBe $ implies $
\xBc X \xfA \xfA Z \xBe,$ nor, conversely,
$ \xBc X \xfA \xfA Z \xBe $ implies $ \xBc X \xfA Y \xfA Z \xBe.$

Consider $I:=\{x,y,z\}.$

(1.1) Let $U:=\{ \xBc 0,0,0 \xBe, \xBc 1,1,1 \xBe, \xBc 0,1,0 \xBe,
\xBc 1,0,1 \xBe, \xBc 1,1,0 \xBe, \xBc 0,0,1 \xBe \}.$
Then $ \xBc x \xfA \xfA z \xBe,$ as all combinations for $x$ and $y$
exist, i.e. paths with
the projections $ \xBc 0,0 \xBe,$ $ \xBc 0,1 \xBe,$ $ \xBc 1,0 \xBe,$ $
\xBc 1,1 \xBe.$
Fix, e.g., $y=1.$ Then the paths through $y=1$ are
$ \xBc 1,1,1 \xBe,$ $ \xBc 0,1,0 \xBe,$ $ \xBc 1,1,0 \xBe,$ but $ \xBc
0,1,1 \xBe $ is missing. So
$ \xBc x \xfA y \xfA z \xBe $ does not hold.

(1.2) Let $U:=\{ \xBc 0,0,0 \xBe, \xBc 1,1,1 \xBe \}.$ Then $ \xBc x \xfA
\xfA z \xBe $ trivially fails, but
$ \xBc x \xfA y \xfA z \xBe $ holds.

(2)

Consider $I:=\{x,a,b,c,d,z\}.$

Let $ \xbS:=\{111111,$ 011110, 011101, 111100, 110111, $010000\}.$

Then $ \xCN x(abcd)z,$ $x(abc)z,$ $ \xCN x(ab)z.$

For $ \xCN x(abcd)z,$ fix $abcd=1111,$ then $111111,011110 \xbe \xbS,$
but, e.g., $011111 \xce \xbS.$

For $x(abc)z,$ the following combinations of abc exist: $111,101,100.$
The result is trivial for 101 and 100. For 111, all combinations for $x$
and $z$
with 0 and 1 exist.

For $ \xCN x(ab)z,$ fix $ab=10,$ then $110111,010000 \xbe \xbS,$ but
there is, e.g., no
$110xy0 \xce \xbS.$

See Diagram \ref{Diagram Change-Rel} (page \pageref{Diagram Change-Rel})

$ \xcz $
\\[3ex]

\ee

$ \xCO $

\vspace{10mm}

\begin{diagram}

\label{Diagram Change-Rel}

\centering
\setlength{\unitlength}{1mm}
{\renewcommand{\dashlinestretch}{30}
\begin{picture}(160,190)(0,0)

\put(10,170){$x$}
\put(30,170){$a$}
\put(50,170){$b$}
\put(70,170){$c$}
\put(90,170){$d$}
\put(110,170){$z$}

\put(10,167){\circle*{1}}
\put(30,167){\circle*{1}}
\put(50,167){\circle*{1}}
\put(70,167){\circle*{1}}
\put(90,167){\circle*{1}}
\put(110,167){\circle*{1}}

\path(10,160)(110,160)
\path(10,140)(30,159)(90,159)(110,140)

\put(120,150){$\xCN  \xBc x \xfA abcd \xfA z \xBe $}
\put(4,150){(1)}

\path(10,110)(70,110)(90,90)(110,90)
\path(10,90)(30,109)(70,109)(90,89)(110,109)

\put(120,100){$  \xBc x \xfA abc \xfA z \xBe $}
\put(4,100){(2)}

\put(10,78){add paths equal on $abc$, different on $d$,
to compensate lacking paths in (1)}

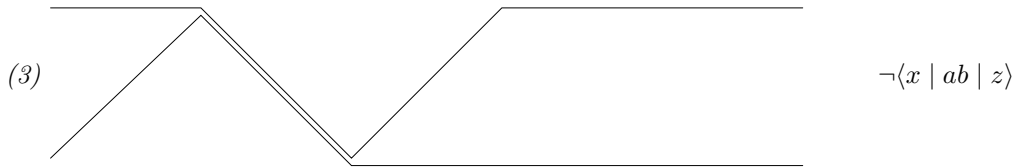
\path(10,60)(30,60)(50,40)(70,60)(110,60)
\path(10,40)(30,59)(50,39)(110,39)

\put(120,50){$\xCN  \xBc x \xfA ab \xfA z \xBe $}
\put(4,50){(3)}

\put(10,28){add paths different on $ab$, singletons on $c$,
so they don't disturb on $abc$:}
\put(10,24){seen on $abc$, the added paths are singletons, so
they respect automatically}
\put(10,20){$  \xBc x \xfA abc \xfA z \xBe $}

\end{picture}
}

\end{diagram}

\vspace{4mm}

$ \xCO $
\subsubsection{
Example of a rule derived from the basic rules
}

We will use the following definition.

\bd

$\hspace{0.01em}$


\label{Definition Func-Sigma-Mu}

Given $ \xbS $ as above, set

$ \xbS_{ \xbm }:=\{ \xBc X,Y,Z \xBe:$ $X,Y,Z$ are pairwise disjoint
subsets of $I,$ $ \xBc X \xfA Y \xfA Z \xBe \xce \xbS,$ but
for all $X' \xcb X$ and all $Z' \xcb Z$ $ \xBc X' \xfA Y \xfA Z \xBe \xbe
\xbS $ and $ \xBc X \xfA Y \xfA Z' \xBe \xbe \xbS \}.$

We will sometimes write $ \xBc X,X' \xfA Y \xfA Z \xBe $ etc. for $ \xBc X
\xcv X' \xfA Y \xfA Z \xBe.$

When we write $ \xBc X,X' \xfA Y \xfA Z \xBe $ etc., we will tacitly
assume that all sets
$X,X',Y,Z$ are pairwise disjoint.

\ed

\br

$\hspace{0.01em}$


\label{Remark Func-Sigma-Mu}

(1)
$ \xbS_{ \xbm }$ contain thus the minimal $X$ and $Z$ for fixed $Y,$ such
that $ \xBc X \xfA Y \xfA Z \xBe \xce \xbS.$

(2)
By rule (b), for all $ \xBc X \xfA Y \xfA Z \xBe \xbe \xbS,$ there is $
\xBc X',Y,Z' \xBe \xbe \xbS_{ \xbm }$ $X \xcc X',$ $Z \xcc Z',$
unless all $ \xbs,$ $ \xbt $ such that $ \xbs \xex Y= \xbt \xex Y$ can be
combined.

\er

As the cases can become a bit complicated, it is important to develop
a good intuition and representation of the problem. We do this now
in the proof of the following fact, where we use the result we want to
prove
to guide our intuition.

\bfa

$\hspace{0.01em}$


\label{Fact Func-Complic}

Let $ \xbS $ be closed under rules $(a)-(d).$
Then, if $ \xBc X,X',X'' \xfA Y \xfA Z,Z',Z'' \xBe \xbe \xbS_{ \xbm },$
then $ \xBc X,Z' \xfA X',Y,Z'' \xfA X'',Z \xBe \xce \xbS.$

\efa

\subparagraph{
Proof
}

$\hspace{0.01em}$


$ \xCO $

\vspace{10mm}

\begin{diagram}

\label{Diagram Func-Compos}

\centering
\setlength{\unitlength}{1mm}
{\renewcommand{\dashlinestretch}{30}
\begin{picture}(160,100)(0,40)

\path(10,120)(115,120)

\path(10,122)(10,118)
\path(25,122)(25,118)
\path(40,122)(40,118)
\path(55,122)(55,118)
\path(70,122)(70,118)
\path(85,122)(85,118)
\path(100,122)(100,118)
\path(115,122)(115,118)

\put(16,125){$X$}
\put(31,125){$X'$}
\put(46,125){$X''$}
\put(61,125){$Y$}
\put(76,125){$Z$}
\put(91,125){$Z'$}
\put(106,125){$Z''$}

\put(15,112){$\xbs_{X}$}
\put(30,112){$\xbs_{X'}$}
\put(45,112){$\xbs_{X''}$}
\put(56,112){$\xbs_{Y}=\xbt_{Y}$}
\put(75,112){$\xbt_{Z}$}
\put(90,112){$\xbt_{Z'}$}
\put(105,112){$\xbt_{Z''}$}

\path(10,110)(115,110)

\path(25,100)(115,100)
\put(118,100){(1)}

\path(10,95)(25,95)
\path(40,95)(115,95)
\put(118,95){(2)}

\path(10,90)(40,90)
\path(55,90)(115,90)
\put(118,90){(3)}

\path(10,85)(70,85)
\path(85,85)(115,85)
\put(118,85){(4)}

\path(10,80)(85,80)
\path(100,80)(115,80)
\put(118,80){(5)}

\path(10,75)(100,75)
\put(118,75){(6)}

\path(25,65)(40,65)
\path(55,65)(70,65)
\path(100,65)(115,65)
\path(25,64)(40,64)
\path(55,64)(70,64)
\path(100,64)(115,64)

\put(30,57)
{Prerequisite:
$\xbs_{X'}=\xbt_{X'}$, $\xbs_{Y}=\xbt_{Y}$, $\xbs_{Z''}=\xbt_{Z''}$}

\end{picture}
}

\end{diagram}

\vspace{4mm}

$ \xCO $

The upper line is the final aim.
Line (1) expresses that we can combine all parts except $s_{X},$ by
$ \xBc X',X'' \xfA Y \xfA Z,Z',Z'' \xBe,$ which holds by $ \xBc X,X'
,X'' \xfA Y \xfA Z,Z',Z'' \xBe \xbe \xbS_{ \xbm },$
by similar arguments, we can combine as indicated in lines
$(2)-(6).$
We now assume $ \xBc X,Z' \xfA X',Y,Z'' \xfA X'',Z \xBe \xbe \xbS.$ So
we have to look at fragments,
which agree on $X',Y,Z''.$ This is, for instance, true for (1) and (3).

We turn this argument now into a formal proof:

Assume

(A) $ \xBc X,Z' \xfA X',Y,Z'' \xfA X'',Z \xBe \xbe \xbS,$ and

(B) $ \xBc X,X',X'' \xfA Y \xfA Z,Z',Z'' \xBe \xbe \xbS_{ \xbm }.$

(C) $ \xBc X,X' \xfA Y \xfA Z,Z',Z'' \xBe $ by (B), see line (3)

(D) $ \xBc X \xfA X',Y,Z',Z'' \xfA X'',Z \xBe $ by (A) and rule (c)

(E) $ \xBc X \xfA X',Y \xfA Z,Z',Z'' \xBe $ by (C) and rule (c)

(F) $ \xBc X \xfA X',Y \xfA Z',Z'' \xBe $ by (E) and (b)

(G) $ \xBc X \xfA X',Y \xfA X'',Z,Z',Z'' \xBe $ by (D) and (F) and (d)

(K) $ \xBc X \xfA X',X'',Y \xfA Z,Z',Z'' \xBe $ by (G) and (c)

(L) $ \xBc X',X'' \xfA Y \xfA Z,Z',Z'' \xBe $ by (B), see line (1)

(M) $ \xBc Z,Z',Z'' \xfA X',X'',Y \xfA X \xBe $ by (K) and (a)

(N) $ \xBc Z,Z',Z'' \xfA Y \xfA X',X'' \xBe $ by (L) and (a)

(O) $ \xBc Z,Z',Z'' \xfA Y \xfA X,X',X'' \xBe $ by (M) and (N) and (d)

(P) $ \xBc X,X',X'' \xfA Y \xfA Z,Z',Z'' \xBe $ by (O) and (a).

So we conclude $ \xBc X,X',X'' \xfA Y \xfA Z,Z',Z'' \xBe \xbe \xbS,$ a
contradiction.

Comment:

We first move $Z',Z'' $ to the right, and then $X',X'' $ to the left.

Moving $Z',Z'':$

We use $X'' $ (or $Z)$ on the right, which not be changed, therefore we
can use
line (3), resulting in

(C) $ \xBc X,X' \xfA Y \xfA Z,Z',Z'' \xBe,$ or, directly

$(C')$ $ \xBc X,X' \xfA Y \xfA Z',Z'' \xBe,$ again by $ \xbS_{ \xbm },$

which is modified to

(F) $ \xBc X \xfA X',Y \xfA Z',Z'' \xBe,$ so we have on the right $Z'
,Z'' $ which we want to move.

We put $Z' $ in the middle ($Z'' $ is there already) of (A), resulting in

(D) $ \xBc X \xfA X',Y,Z',Z'' \xfA X'',Z \xBe.$

Now we can apply (d) to (D) and (F), and have moved $Z',Z'' $ to the
right:

(G) $ \xBc X \xfA X',Y \xfA X'',Z,Z',Z'' \xBe.$

We still have to move $X' $ and $X'' $ to the left of (G), and do this in
an
analogous way.

$ \xcz $
\\[3ex]

Note that our results stays valid, if some of the $X',X'',Z',Z'' $ are
empty.

Aber resultat darf nicht links oder rechts $ \xCQ $ sein.

\bco

$\hspace{0.01em}$


\label{Corollary Func-Complic}

Let $ \xbS $ be closed under rules $(a)-(d).$
Then, if $ \xBc X,X',X'' \xfA Y,Y',Y'' \xfA Z,Z',Z'' \xBe \xbe \xbS_{
\xbm },$ then
$ \xBc X,Y',Z' \xfA X',Y,Z'' \xfA X'',Y'',Z \xBe \xce \xbS.$

Thus, if, for given $Y \xcv Y' \xcv Y'',$
$ \xBc X,X',X'' \xfA Y,Y',Y'' \xfA Z,Z',Z'' \xBe \xbe \xbS_{ \xbm },$
then for no distribution of
$X \xcv X' \xcv X'' \xcv Y \xcv Y' \xcv Y'' \xcv Z \xcv Z' \xcv Z'' $ such
that the outward elements are
non-empty,
$ \xBc X,Y',Z' \xfA X',Y,Z'' \xfA X'',Y'',Z \xBe \xbe \xbS.$

\eco

\subparagraph{
Proof
}

$\hspace{0.01em}$


Suppose $ \xBc X,Y',Z' \xfA X',Y,Z'' \xfA X'',Y'',Z \xBe \xbe \xbS.$
Then by rule (c)
$ \xBc X,Z' \xfA X',Y,Y',Y'',Z'' \xfA X'',Z \xBe \xbe \xbS.$ Set
$Y_{1}:=Y \xcv Y' \xcv Y''.$ Then
$ \xBc X,Z' \xfA X',Y_{1},Z'' \xfA X'',Z \xBe \xbe \xbS,$ and
$ \xBc X,X',X'' \xfA Y_{1} \xfA Z,Z',Z'' \xBe \xbe \xbS_{ \xbm },$
contradicting
Fact \ref{Fact Func-Complic} (page \pageref{Fact Func-Complic}).
$ \xcz $
\\[3ex]
\clearpage
\subsection{
Examples of new rules
}
\subsubsection{
New rules
}

Above rules $(a)-(d)$ are not the only ones to hold, and we introduce now
more complicated ones, and show that they hold in our situation.
Of the possibly infinitary rules, only (Loop1) is given in full
generality, (Loop2) is only given to illustrate
that even the infinitary rule (Loop1) is not all there is.

For warming up, we consider the following short version of (Loop1):

\be

$\hspace{0.01em}$


\label{Example Loop1-Simple}

$ABC,ACD,ADE,AEB \xch ABE.$

We show that this rule holds in all $ \xbS.$

Suppose $A: \xbs,$ $B: \xbs = \xbt,$ $C: \xbt,$ so by $ \xCf ABC,$
there is $ \xbr_{1}$ such that

$A: \xbr_{1}= \xbs,$ $B: \xbr_{1}= \xbs = \xbt,$ $C: \xbr_{1}= \xbt.$
So by $ \xCf ACD,$ there is $ \xbr_{2}$ such that

$A: \xbr_{2}= \xbs,$ $C: \xbr_{2}= \xbr_{1}= \xbt,$ $D: \xbr_{2}= \xbt
.$ So by $ \xCf ADE,$ there is $ \xbr_{3}$ such that

$A: \xbr_{3}= \xbs,$ $D: \xbr_{3}= \xbr_{2}= \xbt,$ $E: \xbr_{3}= \xbt
.$ So by $ \xCf AEB,$ there is $ \xbr_{4}$ such that

$A: \xbr_{4}= \xbs,$ $E: \xbr_{4}= \xbr_{3}= \xbt,$ $B: \xbr_{4}= \xbt =
\xbs.$

So $ \xCf ABE.$

We abbreviate this reasoning by:

(1) $ \xCf ABC:$ $A: \xbs,$ $B: \xbs = \xbt,$ $C: \xbt $

(2) $ \xCf ACD:$ $(1)+ \xbt $

(3) $ \xCf ADE:$ $(2)+ \xbt $

(4) $ \xCf AEB:$ $(3)+ \xbt $

So $ \xCf ABE.$

It is helpful to draw a little diagram as in the following
Table \ref{Table Loop1-Short} (page \pageref{Table Loop1-Short}).

\ee

.
\begin{table}

\label{Table Loop1-Short}
\begin{center}
\begin{tabular}{|c|c|c|c|c|c|c|}
\hline
\multicolumn{7}{|c|}{Validity of $ ABC,ACD,ADE,AEB \xch ABE $} \\
\hline

 \xEH $ \xCf A$ \xEH $B$ \xEH $C$ \xEH $D$ \xEH $E$ \xEH \xEP
 \xEH $ \xbs $ \xEH $ \xbs = \xbt $ \xEH \xEH \xEH $ \xbt $ \xEH $ \xCf
ABE$? \xEP
\hline

(1) $ \xbr_{1}$ \xEH $ \xbs $ \xEH $ \xbs = \xbt $ \xEH $ \xbt $ \xEH \xEH
\xEH $ \xCf ABC$ \xEP
(2) $ \xbr_{2}$ \xEH $ \xbs $ \xEH \xEH $ \xbr_{1}= \xbt $ \xEH $ \xbt $
\xEH \xEH $ \xCf ACD$ \xEP
(3) $ \xbr_{3}$ \xEH $ \xbs $ \xEH \xEH \xEH $ \xbr_{2}= \xbt $ \xEH $
\xbt $ \xEH $ \xCf ADE$ \xEP
(4) $ \xbr_{4}$ \xEH $ \xbs $ \xEH $ \xbs = \xbt $ \xEH \xEH \xEH $
\xbr_{3}= \xbt $ \xEH $ \xCf AEB$ \xEP
\hline
\end{tabular}
\end{center}
\end{table}

We introduce now some new rules.

\bd

$\hspace{0.01em}$


\label{Definition New-Rules}

 \xEI

 \xDH
(Bin1)

$XYZ,XY' Z,Y(XZ)Y' \xch X(YY')Z$

 \xDH
(Bin2)

$XYZ,XZY',Y(XZ)Y' \xch X(YY')Z$

 \xDH
(Loop1)

$AB_{1}B_{2}, \Xl,AB_{i-1}B_{i},AB_{i}B_{i+1},AB_{i+1}B_{i+2}, \Xl
,AB_{n-1}B_{n},AB_{n}B_{1} \xch AB_{1}B_{n}$
$ \xDB $
so we turn $AB_{n}B_{1}$ around to $AB_{1}B_{n}.$

When we have to be more precise, we will denote this condition
$(Loop1_{n})$ to fix the length.

 \xDH
(Loop2)

$ABC,ACD,DAE,DEF,FDG,FGH,HFB \xch HBF:$

 \xEJ

\ed

The complicated structure of these rules suggests already that the ternary
relations are not the right level of abstraction to speak about
construction
of functions from fragments. This is made formal by our main result below,
which shows that there is no finite characterization by such relations.
In other words, the main things happen behind the screen.

\bfa

$\hspace{0.01em}$


\label{Fact New-Valid}

The new rules are valid in our situation.

\efa

\subparagraph{
Proof
}

$\hspace{0.01em}$


 \xEI

 \xDH
(Bin1)

(1) $ \xCf XYZ:$ $X: \xbs,$ $Y: \xbs = \xbt,$ $Z: \xbt $

(2) $XY' Z:$ $X: \xbs,$ $Y': \xbs = \xbt,$ $Z: \xbt $

(3) $Y(XZ)Y':$ $(1)+(2)$

So $X(YY')Z.$

 \xDH
(Bin2)

Let $X: \xbs,$ $Y: \xbs = \xbt,$ $Y': \xbs = \xbt,$ $Z: \xbt $

(1) $ \xCf XYZ:$ $X: \xbs,$ $Y: \xbs = \xbt,$ $Z: \xbt $

(2) $ \xCf XZY':$ $(1)+ \xbt $

(3) $Y(XZ)Y':$ $(1)+(2)$

So $X(YY')Z.$

 \xDH
(Loop1)

(1) $AB_{1}B_{2}:$ $A: \xbs,$ $B_{1}: \xbs = \xbt,$ $B_{2}: \xbt $

(2) $AB_{2}B_{3}:$ $(1)+ \xbt $

 \Xl.

(i-1) $AB_{i-1}B_{i}:$ $(i-2)+ \xbt $

(i) $AB_{i}B_{i+1}:$ $(i-1)+ \xbt $

$(i+1)$ $AB_{i+1}B_{i+2}:$ $(i)+ \xbt $

 \Xl.

$(\xCf n-1)$ $AB_{n-1}B_{n}:$ $(n-2)+ \xbt $

$(\xCf n)$ $AB_{n}B_{1}:$ $(n-1)+ \xbt $

So $AB_{1}B_{n}.$
 \xDH
(Loop2)

Let

(1) $ \xCf ABC:$ $A: \xbs,$ $B: \xbs = \xbt,$ $C: \xbt $

(2) $ \xCf ACD:$ $1+ \xbt $

(3) $ \xCf DAE:$ $2+ \xbs $

(4) $ \xCf DEF:$ $3+ \xbs $

(5) $ \xCf FDG:$ $4+ \xbt $

(6) $ \xCf FGH:$ $5+ \xbt $

(7) $ \xCf HFB:$ $6+ \xbs $

So $ \xCf HBF$ by $B: \xbs = \xbt.$

 \xEJ

Note that we use here $B: \xbs = \xbt,$ $E: \xbs = \xbt,$ $H: \xbs =
\xbt,$ whereas the other tripels
are used for other functions.

$ \xcz $
\\[3ex]

Next we show that the full (Loop1) cannot be derived from the basic
rules $(a)-(d)$ and (Bin1), and shorter versions of (Loop1).
(This is also a consequence of the sequel, but we want to point it out
right away.)

\bfa

$\hspace{0.01em}$


\label{Fact Loop-Indep}

Let $n \xcg 1,$ then $(Loop1_{n})$ does not follow from the rules
$(a)-(d),$ $(\xCQ),$ (Bin1), and the shorter versions of (Loop1)

\efa

\subparagraph{
Proof
}

$\hspace{0.01em}$


Consider the following set of tripels $L \xcv L' $ over
$I:=\{a,b_{1}, \Xl,b_{n}\}:$

$L:=\{ab_{1}b_{2},$  \Xl, $ab_{i}b_{i+1},$  \Xl, $ab_{n-1}b_{n},$
$ab_{n}b_{1}\},$

$L':=\{ \xCQ AB:$ $A \xcs B= \xCQ,$ $A \xcv B \xcc I\},$

and close this set under symmetry (rule (a)).
Call the resulting set $ \xda.$

Note that, on the outside, we have $ \xCQ $ or singletons, inside
singletons or $ \xCQ.$
If the inside is $ \xCQ,$ one of the outside sets must also be $ \xCQ.$

When we look at $L,$ and define a relation $<$ by $x<y$ iff $axy \xbe L,$
we see
that the only $<$-loop is $b_{1}<b_{2}< \Xl <b_{n}<b_{1}.$

We show first that $ \xda $ is closed under rules $(a)-(d)$
(see Definition 
\ref{Definition Basic-Rules} (page 
\pageref{Definition Basic-Rules})).

(a) is trivial.

(b) If $W= \xCQ $ or $Z= \xCQ,$ this is trivial, if $W=Z,$ this is
trivial, too.

(c) If $Z \xcv W= \xCQ,$ this is trivial, if $Z \xcv W$ is a singleton,
so
$Z= \xCQ $ or $W= \xCQ $ or $Z=W.$ $Z= \xCQ $ or $W= \xCQ $ are
trivial, otherwise $Z=W$ contradicts disjointness.

(d) $Z= \xCQ $ is trivial, so is $W= \xCQ,$ otherwise $Z=W$ contradicts
disjointness.

(Bin1)
$X= \xCQ $ or $Z= \xCQ $ are trivial, otherwise $X=Z$ is excluded by
disjointness.
So we are in $L' $ for $Y(XZ)Y'.$ So $Y= \xCQ $ or $Y' = \xCQ $ and it is
trivial.

Obviously, $(Loop1_{n})$ does not hold.

We show now that all $(Loop1_{k}),$ $0 \xck k<n$ hold.

The cases $n=1,$ $n=2$ are trivial.

Consider the case $2<k<n.$

This has the form $AB_{1}B_{2},AB_{2}B_{3}, \Xl,AB_{k-1}B_{k},AB_{k}B_{1}
\xch AB_{1}B_{k}.$

If $A= \xCQ $ or $B_{k}= \xCQ,$ the condition holds.

So assume $A,B_{k} \xEd \xCQ.$ Thus, by above remark, descending to
$B_{k-1}$ etc., we see that
all $B_{i} \xEd \xCQ,$ $1 \xck i \xck k.$
Thus, all prerequisites are in $L.$
Moreover, $ \xCf A$ has to be $ \xCf a,$ which
is the only element occuring repeatedly on the outside. Consider now the
relation $<' $ defined by $U<' V$ iff $ \xCf AUV$ is among the
prerequisites.
We then have $B_{1}<' B_{2}<'  \Xl <' B_{k}<' B_{1},$ where all $B_{i}$
are some $b_{j},$
we see that the resulting $<' $-loop is too short, so the prerequisites
cannot hold, and we have a contradiction.

$ \xcz $
\\[3ex]
\subsection{
There is no finite characterization
}

We turn to our main result.
\subsubsection{
Discussion
}

Consider the following simple, short, loop for illustration:

$ABC,ACD,ADE,AEF,AFG,AGB \xch ABG$ - so we can turn $ \xCf AGB$ around to
$ \xCf ABG.$

Of course, this construction may be arbitrarily long.

The idea is now to make $ \xCf ABG$ false, and, to make it coherent, to
make one
of the interior conditions false, too, say $ \xCf ADE.$ We describe this
situation fully, i.e. enumerate all conditions which hold in such a
situation.
If we make now $ \xCf ADE$ true again, we know this is not valid, so any
(finite)
characterization must say ``NO'' to this. But as it is finite, it cannot
describe all the interior tripels of the type $ \xCf ADE$ in a
sufficiently long loop,
so we just change one of
them which it does not ``see'' to FALSE, and it must give the same answer
NO, so
this fails.

Basically, we cannot describe parts of the loop, as the $< \xfA \xfA
>$-language
is not rich enough to express it, we see only the final outcome.

The problem is to fully describe the situation.
\subsubsection{
Composition of layers
}

A very helpful fact is the following:

\bd

$\hspace{0.01em}$


\label{Definition Sigma-S-H}

Let $ \xbS_{j}$ be function sets over $I$ into some set $K,$ $j \xbe J.$

Let $ \xbS $ $:=$ $\{$ $f:I \xcp K^{J}:$ $f(i)=\{ \xBc f_{j}(i),j \xBe:j
\xbe J,f_{j} \xbe \xbS_{j}\}$ $\}.$

So any $f \xbe \xbS $ has the form
$f(i)= \xBc f_{1}(i),f_{2}(i), \Xl,f_{n}(i) \xBe,$ $f_{m} \xbe \xbS_{m}$
(we may assume $J$ to be finite).

Thus, given $f \xbe \xbS,$ $f_{m} \xbe \xbS_{m}$ is defined.

\ed

\bfa

$\hspace{0.01em}$


\label{Fact Sigma-S-H}

For the above $ \xbS $ $ \xBc A \xfA B \xfA C \xBe $ holds iff it holds for all
$
\xbS_{j}.$

Thus, we can destroy the $ \xBc A \xfA B \xfA C \xBe $ independently, and
collect
the results.

\efa

\subparagraph{
Proof
}

$\hspace{0.01em}$


The proof is trivial, and a direct consequence of the fact that
$f=f' $ iff for all components $f_{j}=f'_{j}.$

Suppose for some $ \xbS_{k},$ $k \xbe J,$ $ \xCN  \xBc A \xfA B \xfA C \xBe.$

So for this $k$ there are $f_{k},f'_{k} \xbe \xbS_{k}$ such that
$f_{k}(B)=f'_{k}(B),$ but
there is no $f''_{k} \xbe \xbS_{k}$ such that $f''_{k}(A)=f_{k}(A),$
$f''_{k}(B)=f_{k}(B)=f'_{k}(B),$ $f''_{k}(C)=f'_{k}(C)$
(or conversely).
Consider now some $h \xbe \xbS $ such that $h_{k}=f_{k},$ and $h' $ is
like $h,$ but $h'_{k}=f'_{k},$ so
also $h' \xbe \xbS.$
Then $h(B)=h' (B),$ but there is no $h'' \xbe \xbS $ such that
$h'' (A)=h(A),$ $h'' (B)=h(B)=h' (B),$ $h'' (C)=h' (C).$

Conversely, suppose $ \xBc A \xfA B \xfA C \xBe $ for all $ \xbS_{j}.$
Let $h,h' \xbe \xbS $ such that $h(B)=h' (B),$ so for all $j \xbe J$
$h_{j}(B)=h'_{j}(B),$ where
$h_{j} \xbe \xbS_{j},$ $h'_{j} \xbe \xbS_{j},$ so there are $h''_{j} \xbe
\xbS_{j}$ with
$h''_{j}(A)=h_{j}(A),$ $h''_{j}(B)=h_{j}(B)=h'_{j}(B),$
$h''_{j}(C)=h'_{j}(C)$ for all $j \xbe J.$
Thus, $h'' $ composed of the $h''_{j}$ is in $ \xbS,$ and $h'' (A)=h(A),$
$h'' (B)=h(B)=h' (B),$ $h'' (C)=h' (C).$

$ \xcz $
\\[3ex]
\subsubsection{
Systematic construction
}

Recall the general form of (Loop1) for singletons:

$ab_{1}b_{2}, \Xl,ab_{i-1}b_{i},ab_{i}b_{i+1},ab_{i+1}b_{i+2}, \Xl
,ab_{n-1}b_{n},ab_{n}b_{1} \xch ab_{1}b_{n}$

We will fully describe a model of above tripels, with the exception of
$ab_{1}b_{n}$
and $ab_{i}b_{i+1}$ which will be made to fail, and all other $ \xBc X
\xfA Y \xfA Z \xBe $
which are not in above list of tripels to preserve, will fail, too
(except for $X= \xCQ $ or $Z= \xCQ).$

Thus, the tripels to preserve are:

$P$ $:=$ $\{ab_{1}b_{2}, \Xl,ab_{i-1}b_{i},$ (BUT NOT $ab_{i}b_{i+1})$
$,ab_{i+1}b_{i+2}, \Xl,ab_{n-1}b_{n},ab_{n}b_{1}\}$

We use the following fact:

\bfa

$\hspace{0.01em}$


\label{Fact Const-All-2-H}

Let $X \xcc I,$ $card(X)>1,$ $ \xbS_{X}$ $:=$ $\{$ $ \xbs:I \xcp \{0,1\}$
: $card\{x \xbe X: \xbs (x)=0\}$ is even $\}$

Then $ \xCN ABC$ iff $A \xcs X \xEd \xCQ,$ $C \xcs X \xEd \xCQ,$ $X \xcc
A \xcv B \xcv C.$

\efa

\subparagraph{
Proof
}

$\hspace{0.01em}$


``$ \xci $'':

Suppose $A \xcs X \xEd \xCQ,$ $C \xcs X \xEd \xCQ,$ $X \xcc A \xcv B
\xcv C.$

Take $ \xbs $ such that $card\{x \xbe X: \xbs (x)=0\}$ is odd, then $ \xbs
\xce \xbS_{X}.$
As $X \xcC A \xcv B,$ there is $ \xbt \xbe \xbS_{X}$ such that $ \xbs \xex
A \xcv B= \xbt \xex A \xcv B.$
As $X \xcC B \xcv C,$ there is $ \xbr \xbe \xbS_{X}$ such that $ \xbr \xex
B \xcv C= \xbs \xex B \xcv C.$
Thus, $ \xbt \xex B= \xbr \xex B.$
If there were $ \xba \xbe \xbS_{X}$ such that $ \xba \xex A \xcv B= \xbt
\xex A \xcv B$ and $ \xba \xex B \xcv C= \xbr \xex B \xcv C,$ then
$ \xba \xex A \xcv B \xcv C= \xbs \xex A \xcv B \xcv C,$ contradiction

``$ \xch $'':

Suppose $A \xcs X= \xCQ $ or $C \xcs X= \xCQ,$ or $X \xcC A \xcv B \xcv
C.$ We show $ \xCf ABC.$

Case 1:
$C \xcs X= \xCQ.$ Let $ \xbs, \xbt \xbe \xbS_{X}$ such that $ \xbs \xex
B= \xbt \xex B.$
As $C \xcs X= \xCQ,$ we can continue $ \xbs \xex A \xcv B$ as we like.

Case 2, $A \xcs X= \xCQ,$ analogous.

Case 3:
$X \xcC A \xcv B \xcv C.$ But then there is no restriction in $A \xcv B
\xcv C.$

$ \xcz $
\\[3ex]

We will have to make $ab_{1}b_{n}$ false, but $ab_{n}b_{1}$ true. On the
other hand,
we will make $ab_{1}b_{3}$ false, but $ab_{3}b_{1}$ need not be preserved.

This leads to the following definition, which helps to put order into
the cases.

\bd

$\hspace{0.01em}$


\label{Definition dmin-H}

Suppose we have to destroy $ \xCf axy.$ Then

$dmin(axy)$ $:=$ $min\{d(\{a,x,y\},\{a,u,v\}):$ $ \xCf auv$ has to be
preserved $\}$ - $d$ the
counting Hamming distance.

\ed

Thus, $dmin(ab_{1}b_{n})=0$ (as $ab_{n}b_{1}$ has to be preserved),
$dmin(ab_{1}b_{3})=1$
(because $ab_{1}b_{2}$ has to be preserved, but not $ab_{3}b_{1}).$

We introduce the following order defined from the loop prerequisites to be
preserved.

\bd

$\hspace{0.01em}$


\label{Definition Loop-Order-H}

Order the elements by following the string of sequences to be preserved as
follows:

$b_{i+1} \xeb b_{i+2} \xeb  \Xl  \xeb b_{n-1} \xeb b_{n} \xeb b_{1} \xeb
b_{2} \xeb  \Xl  \xeb b_{i-1} \xeb b_{i}$

Note that the interruption at $ab_{i}b_{i+1}$ is crucial here - otherwise,
there
would be a cycle.

As usual, $ \xec $ will stand for $ \xeb $ or $=.$
\subsubsection{
The cases to consider
}

\ed

The elements to consider are: $a,b_{1}, \Xl,b_{n}.$

Recall that the tripels to preserve are:

$P$ $:=$ $\{ab_{1}b_{2}, \Xl,ab_{i-1}b_{i},$ (BUT NOT $ab_{i}b_{i+1})$
$,ab_{i+1}b_{i+2}, \Xl,ab_{n-1}b_{n},ab_{n}b_{1}\}$

The $ \xBc X \xfA Y \xfA Z \xBe $ to destroy are (except when $X= \xCQ $
or $Z= \xCQ):$

 \xEh

 \xDH
all $ \xBc X \xfA \xfA Z \xBe $

 \xDH
all $ \xBc X \xfA Y \xfA Z \xBe $ such that $X \xcv Y \xcv Z$ has $>3$
elements

 \xDH
all tripels which do not have $ \xCf a$ on the outside, e.g. $ \xCf bgc$

 \xDH
and the following tripels:

(the (0) will be explained below - for the moment, just ignore it)

$ab_{1}b_{3}, \Xl,ab_{1}b_{n-1},$ $ab_{1}b_{n}$ (0)

$ab_{2}b_{1}$ (0), $ab_{2}b_{4}, \Xl,ab_{2}b_{n}$

$ab_{3}b_{1},$ $ab_{3}b_{2}$ (0), $ab_{3}b_{5}, \Xl,ab_{3}b_{n}$

 \Xl.

$ab_{i}b_{1},$ $ab_{i}b_{2},, \Xl,$ ALSO $ab_{i}b_{i+1}, \Xl
,ab_{i}b_{n}$

 \Xl.

$ab_{n-2}b_{1},, \Xl,$ $ab_{n-2}b_{n-3}$ (0), $ab_{n-2}b_{n}$

$ab_{n-1}b_{1},, \Xl,$ $ab_{n-1}b_{n-2}$ (0),

$ab_{n}b_{1},$ $, \Xl,ab_{n}b_{n-1}$ (0)

 \xEj
\subsubsection{
Solution of the cases
}

We show how to destroy all tripels mentioned above, while preserving
all tripels in $P.$

 \xEh
 \xDH
all $ \xBc X \xfA Y \xfA Z \xBe $ where $X \xcv Y \xcv Z$ has $>3$
elements:

See Fact 
\ref{Fact Const-All-2-H} (page 
\pageref{Fact Const-All-2-H})  with the $X$ there with 4 elements,
for all such $X,Y,Z$
separately, so all tripels in $P$ are preserved.

 \xDH
all $ \xBc X \xfA Y \xfA Z \xBe $ with 1 element: -

 \xDH
all $ \xBc X \xfA \xfA Z \xBe:$

This can be
done by considering $ \xbS_{j}:=\{0_{c},1_{c}\}.$ Then, say for
$a,c,$ we have to examine the fragments 00 and 11, but there is no 10 or
01.
For $ \xBc a \xfA b \xfA c \xBe $ this is no problem, as we have only the two
000,
111, which
do not agree on $b.$

 \xDH
all $ \xBc X \xfA Y \xfA Z \xBe $ with 2 elements: eliminated by $ \xBc X
\xfA \xfA Z \xBe $

 \xDH
all $ \xBc X \xfA Y \xfA Z \xBe $ with 3 elements:

 \xEh
 \xDH
$ \xCf a$ is not on the outside
 \xEh
 \xDH
$ \xCf a$ is in the middle, we need
$ \xCN xay:$ Consider $ \xbS $ with 2 functions, $0_{c},$ and the second
defined by
$a=0,$ and all $u=1$ for $u \xEd a.$
Obviously, $ \xCN xay.$ Recall that all tripels to be preserved have $
\xCf a$ on the
outside, and some other element $x$ in the middle. Then the two functions
are
different on $x.$

 \xDH
$ \xCf a$ is not in $ \xCf xyz,$ we need $ \xCN xyz:$
Consider $ \xbS $ with 2 functions, $0_{c},$ and the second defined by
$a=y=0,$ all $u=1$
for $u \xEd a,$ $u \xEd y.$
As $ \xCf a$ is neither $x$ nor $z,$ $ \xCN xyz.$ If some $ \xCf uvw$ has
$ \xCf a$ on the outside, say
$u=a,$ then both functions are 000 or 0vw on this tripel, so $ \xCf uvw$
holds.

 \xEj

 \xDH

$ \xCf a$ is on the outside, we destroy $ \xCf ayz:$

 \xEh

 \xDH
Case $dmin(ayz)>0$:

Take as $ \xbS $ the set of all functions with values in $\{0,1\},$ but
eliminate
those with $a=y=z=0.$
Then $ \xCN ayz$ (we have $100,001,101,$ but not 000),
but for all $ \xCf auv$ with $d(\{a,y,z\},\{a,u,v\})>0$ $ \xCf auv$ has
all possible
combinations,
as all combinations for $ \xCf ay$ and $ \xCf az$ exist.

 \xDH
Case $dmin(ayz)=0.$

The elements with $dmin=0$ are:

$ab_{1}b_{n},$ $ab_{2}b_{1},$  \Xl, $ab_{i}b_{i-1},$ NOT
$ab_{i+1}b_{i},$ $ab_{i+2}b_{i+1},$  \Xl, $ab_{n-1}b_{n-2},$
$ab_{n}b_{n-1},$
they were marked with (0) above.

$ \xbS $ will again have 2 functions, the first is always $0_{c}.$

The second function: Always set $a=1.$

We see that the tripels with $dmin=0$ to be destroyed have the form $ \xCf
ayz,$ where
$z$ is the immediate $ \xeb $-predecessor of $y$ in
above order - see Definition 
\ref{Definition Loop-Order-H} (page 
\pageref{Definition Loop-Order-H}).
Conversely, those to be preserved (in $P)$ have the form $ \xCf azy,$
where again
$z$ is the immediate $ \xeb $-predecessor of $y.$

We set $z' =1$ for all $z' \xec z,$ and $y' =0$ for all $y' \xed y.$
Recall that
$z \xeb y,$ so we have the picture
$b_{i+1}=1, \Xl,z=1,y=0, \Xl,b_{i}=0.$

Then $ \xCN ayz,$ as we have the fragments 000, 101. But $ \xCf azy,$ as
we have the
fragments 000, 110.
Moreover, considering the successors of the sequence, we give the values
11, or
10, or 00. This results in the function fragments for $ \xCf auv$ as 111,
or 110, or
100. But the resulting fragment sets (together with $0_{c})$ are then:
$\{000,111\},$ $\{000,110\},$ $\{000,100\}.$ They all make $ \xCf auv$
true. Thus, all
tripels in $P$ are preserved.

 \xEj
 \xEj
 \xEj
\clearpage
\subsection{
Systematic construction of new rules
}

This section is an outline - not a formal proof - for constructing a
complete rule set for our scenario.

We give here a general way how to construct new rules of the type
ABC, DEF,  \Xl. $ \xch $ XYZ which are valid in our situation.
\subsubsection{
Consequences of a single tripel
}

Let $(XX' X'')Y(ZZ' Z'')$ be a tripel, then all consequences of this
single tripel
have the form $X(X' YZ')Z$ (up to symmetry).

Obviously, such $X(X' YZ')Z$ are consequences, using rules (b) and (c).

We now give counterexamples to other forms, to show that they are not
consequences in our setting.
We always assume that the outside is not $ \xCQ.$
We consider $A=B=C=\{0,1\},$ and subsets of $A \xCK B \xCK C.$

 \xEh
 \xDH
$Y$ decreases:

Consider $\{000,111\},$ then ABC, but not $A \xCQ C.$

 \xDH
$Z$ increases:

Consider $\{000,101\},$ then $A \xCQ B,$ but not $A \xCQ (BC).$

 \xDH
$X$ goes from left to right:

Consider $\{000,110\},$ then (AB)C, but not $A(BC)$

 \xDH
$Y$ increases by some arbitrary $W:$

Consider $\{000,101,110,011\},$ then $A \xCQ C,$ but not ABC.

 \xEj
\subsubsection{
Construction of function trees
}

We can construct new functions from two old functions using tripels ABC,
so,
in a more general way, we have a binary function construction tree, where
the old functions are the leaves, and the new function is the root.
The form of such a tree is obvious, the tripels used are either
directly given, or consequences of such tripels. In
Example \ref{Example R-3} (page \pageref{Example R-3}), for
instance, in the construction of $ \xbr_{2},$ we used ACD, but we could
also
have used e.g. $AC(DD'),$ for some $D'.$
\subsubsection{
Derivation trees
}

Not all such function construction trees are proof trees for a rule
$T_{1}, \Xl,T_{n} \xch T,$ where the $T_{i}$ and $T$ are tripels.

We have to look at the logical structure of the tripels to see what we
need.
In order to show $T=ABC,$ we assume given two arbitrary functions $ \xbs $
and $ \xbt,$
which agree on $B,$ and construct $ \xbr $ such that on A $ \xbr = \xbs,$
on $B$ $ \xbr = \xbs = \xbt $
(the latter, $ \xbs = \xbt $ by prerequisite), and on $C$ $ \xbr = \xbt.$
We will write this as
$A: \xbr = \xbs,$ $B: \xbr = \xbs = \xbt,$ $C: \xbr = \xbt.$

Thus, we have no functions at the beginning, except $ \xbs $ and $ \xbt,$
so all leaves
in a proof tree for $T_{1}, \Xl,T_{n} \xch T$ have to be $ \xbs $ or $
\xbt.$ Moreoever, all we know
about $ \xbs $ and $ \xbt $ is that they agree on $B.$ Thus, we can only
use some
$T_{i}' =A' B' C' $ on $ \xbs $ and $ \xbt $ if $B' \xcc B.$ Likewise, in
the interior of the tree,
we can only use $ \xbs \xex B= \xbt \xex B,$ and, of course, all
equalities which hold be
construction. E.g., in
Example 
\ref{Example R-3} (page 
\pageref{Example R-3}), in the construction of $ \xbr_{2},$ by
construction of $ \xbr_{1},$ $C: \xbr_{1}= \xbt,$ so we can use ACD to
construct $ \xbr_{2}$ from
$ \xbr_{1}$ and $ \xbt.$

At the root, we must have a function $ \xbr $ of the form
$A: \xbr = \xbs,$ $B: \xbr = \xbs = \xbt,$ $C: \xbr = \xbt.$
In Example \ref{Example R-3} (page \pageref{Example R-3}),
$ \xbr_{4},$ at the root, was constructed using AEB from $ \xbr_{3}$ and $
\xbt.$
But we do not interpret $ \xbr_{4}$ as AEB, but as ABE, which is possible,
as $A: \xbr_{4}= \xbs,$ $B: \xbr_{4}= \xbs = \xbt,$ $E: \xbr_{4}= \xbt
.$

Intermediate nodes can be read as an intermediate result $A' B' C' $ by
the same
criteria: They must be functions $ \xbr' $ such that
$A': \xbr' = \xbs,$ $B': \xbr' = \xbs = \xbt,$ $C': \xbr' = \xbt $
and all $B'' $ such that
$B'': \xbs = \xbt $ used up to this node must be subsets of $B',$ as $B'
: \xbs = \xbt $ is the
only hypothesis we then have.
\clearpage
\subsubsection{
Examples
}

$ \xCO $

\vspace{10mm}

\begin{diagram}

\label{Diagram Rule-Const-1}

\centering
\setlength{\unitlength}{1mm}
{\renewcommand{\dashlinestretch}{30}
\begin{picture}(160,190)(0,70)

\put(60,182){Example \ref{Example R-1}}

\put(8,182){$\xbs$}
\put(42,182){$\xbt$}
\path(10,180)(25,165)
\path(27,165)(42,180)
\put(25,162){$\xbr_1$}

\put(60,162){$\xbt$}
\path(28,160)(43,145)
\path(45,145)(60,160)
\put(43,142){$\xbr_2$}

\put(90,122){Examples \ref{Example R-2} and \ref{Example R-4}}

\put(8,122){$\xbs$}
\put(42,122){$\xbt$}
\path(10,120)(25,105)
\path(27,105)(42,120)
\put(25,102){$\xbr_1$}

\put(53,122){$\xbs$}
\put(87,122){$\xbt$}
\path(55,120)(70,105)
\path(72,105)(87,120)
\put(69,102){$\xbr_2$}

\path(28,100)(48,80)
\path(50,80)(70,100)
\put(48,77){$\xbr_3$}

\end{picture}
}

\end{diagram}

\vspace{4mm}

$ \xCO $

$ \xCO $

\vspace{10mm}

\begin{diagram}

\label{Diagram Rule-Const-2}

\centering
\setlength{\unitlength}{1mm}
{\renewcommand{\dashlinestretch}{30}
\begin{picture}(160,195)(0,85)

\put(30,187){Example \ref{Example R-3}}

\put(8,182){$\xbs$}
\put(42,182){$\xbt$}
\path(10,180)(25,165)
\path(27,165)(42,180)
\put(25,162){$\xbr_1$}
\put(50,162){$\xbr_1$, using ABC - $A:\xbs, B:\xbs=\xbt, C:\xbt$}

\put(43,162){$\xbt$}
\path(26,160)(26,145)
\path(28,145)(41,160)
\put(25,142){$\xbr_2$}
\put(50,142){$\xbr_2$, using ACD - $A:\xbr_1=\xbs, C:\xbr_1=\xbt, D:\xbt$}

\put(43,142){$\xbt$}
\path(26,140)(26,125)
\path(28,125)(41,140)
\put(25,122){$\xbr_3$}
\put(50,122){$\xbr_3$, using ADE - $A:\xbr_2=\xbs, D:\xbr_2=\xbt, E:\xbt$}

\put(43,122){$\xbt$}
\path(26,120)(26,105)
\path(28,105)(41,120)
\put(25,102){$\xbr_4$}
\put(50,102){$\xbr_4$, using AEB - $A:\xbr_3=\xbs, E:\xbr_3=\xbt, B:\xbs=\xbt$}
\put(50,95){Interpretation: $ABE$, common part $B:\xbs=\xbt$}

\end{picture}
}

\end{diagram}

\vspace{4mm}

$ \xCO $

$ \xCO $

\vspace{10mm}

\begin{diagram}

\label{Diagram Rule-Const-3}

\centering
\setlength{\unitlength}{0.5mm}
{\renewcommand{\dashlinestretch}{30}
\begin{picture}(320,380)(0,200)

\put(90,372){Example \ref{Example R-5}}

\put(8,362){$\xbs$}
\put(42,362){$\xbt$}
\path(10,360)(25,345)
\path(27,345)(42,360)
\put(24,341){$\xbr_1$}

\put(53,362){$\xbs$}
\put(87,362){$\xbt$}
\path(55,360)(70,345)
\path(72,345)(87,360)
\put(68,341){$\xbr_2$}

\path(28,339)(48,319)
\path(50,319)(70,339)
\put(47,315){$\xbr_3$}

\put(108,362){$\xbs$}
\put(142,362){$\xbt$}
\path(110,360)(125,345)
\path(127,345)(142,360)
\put(124,341){$\xbr_1$}

\put(153,362){$\xbs$}
\put(187,362){$\xbt$}
\path(155,360)(170,345)
\path(172,345)(187,360)
\put(168,341){$\xbr_2$}

\path(128,339)(148,319)
\path(150,319)(170,339)
\put(147,315){$\xbr_4$}

\path(51,312)(96,267)
\path(100,267)(147,312)
\put(97,261){$\xbr_5$}

\put(153,282){$\xbs$}
\put(187,282){$\xbt$}
\path(155,280)(170,265)
\path(172,265)(187,280)
\put(168,261){$\xbr_1$}

\path(103,258)(135,227)
\path(142,227)(166,258)
\put(136,222){$\xbr_6$}

\end{picture}
}

\end{diagram}

\vspace{4mm}

$ \xCO $
\clearpage

Explanation:

By ``prerequisite'' of $ \xbr_{i}$ we mean the set $X$ we used in the
construction,
where $X: \xbs = \xbt.$ For instance, in the construction of $ \xbr_{2}$
in
Example \ref{Example R-1} (page \pageref{Example R-1}),
we used only that $B \xcv C: \xbr_{1}= \xbt $ by the construction of $
\xbr_{1},$ no additional
use of some $ \xbs = \xbt $ was made.

By ``common part'' of $ \xbr_{i}$ we mean the set $X$ such that $X:
\xbr_{i}= \xbs = \xbt.$

\be

$\hspace{0.01em}$


\label{Example R-1}

(Contraction), ABC, $A(BC)D$ $ \xcp $ $AB(CD)$:

(See Diagram 
\ref{Diagram Rule-Const-1} (page 
\pageref{Diagram Rule-Const-1})  upper part.)

 \xEI
 \xDH
$ \xbr_{1}:$ $A: \xbs,$ $B: \xbs = \xbt,$ $C: \xbt $

generated by $ \xCf ABC$ from $ \xbs,$ $ \xbt $

prerequisite $B,$

common part: $B$

$ \xbr_{1}$ can be interpreted as the (trivial) derived tripel $ \xCf ABC$
 \xDH
$ \xbr_{2}:$ $A: \xbr_{1}= \xbs,$ $B: \xbr_{1}= \xbs = \xbt,$ $C:
\xbr_{1}= \xbt,$ $D: \xbt $

generated by $A(BC)D$ from $ \xbr_{1},$ $ \xbt $

prerequisite -,

common part: $B.$

$ \xbr_{2}$ can be interpreted as a derived tripel by $AB(CD).$

$ \xbr_{2}$ can also be interpreted as a derived tripel by $A(BC)D$ or
$A(BD)C.$
Note that these possibilities can be derived from $AB(CD)$ by rule (c),
Weak Union.
 \xEJ

\ee

\be

$\hspace{0.01em}$


\label{Example R-2}

(Bin1), XYZ, $XY' Z,$ $Y(XZ)Y' $ $ \xch $ $X(YY')Z$:

(See Diagram 
\ref{Diagram Rule-Const-1} (page 
\pageref{Diagram Rule-Const-1})  lower part.)

\ee

 \xEI
 \xDH
$ \xbr_{1}:$ $X: \xbs,$ $Y: \xbs = \xbt,$ $Z: \xbt $

generated by $ \xCf XYZ$ from $ \xbs,$ $ \xbt $

prerequisite $Y$

common part: $Y$
 \xDH
$ \xbr_{2}:$ $X: \xbs,$ $Y': \xbs = \xbt,$ $Z: \xbt $

generated by $XY' Z$ from $ \xbs,$ $ \xbt $

prerequisite $Y' $

common part: $Y' $
 \xDH
$ \xbr_{3}:$ $Y: \xbr_{1}= \xbs = \xbt,$ $X: \xbr_{1}= \xbr_{2}= \xbs,$
$Z: \xbr_{1}= \xbr_{2}= \xbt,$ $Y': \xbr_{2}= \xbs = \xbt $

generated by $Y(XZ)Y' $ from $ \xbr_{1},$ $ \xbr_{2}$

prerequisites -

common part: $ \xCf YY' $

$ \xbr_{3}$ can be interpreted as a derived tripel by $X(YY')Z.$
 \xEJ

\be

$\hspace{0.01em}$


\label{Example R-3}

(Loop1) ABC, ACD, ADE, AEB $ \xch $ ABE:

(See Diagram \ref{Diagram Rule-Const-2} (page \pageref{Diagram Rule-Const-2}).)

 \xEI
 \xDH
$ \xbr_{1}:$ $A: \xbs,$ $B: \xbs = \xbt,$ $C: \xbt $

generated by $ \xCf ABC$ from $ \xbs,$ $ \xbt $

prerequisite $B$

common part $B$
 \xDH
$ \xbr_{2}:$ $A: \xbr_{1}= \xbs,$ $C: \xbr_{1}= \xbt,$ $D: \xbt $

generated by $ \xCf ACD$ from $ \xbr_{1},$ $ \xbt $

prerequisite -

common part -

$ \xbr_{2}$ cannot be interpreted as a derived tripel, as there was a
prerequisite
used in its derivation (B), but the common part in $ \xbr_{2}$ is $ \xCQ
.$
 \xDH
$ \xbr_{3}$ similar to $ \xbr_{2}:$

$ \xbr_{3}:$ $A: \xbr_{2}= \xbs,$ $D: \xbr_{2}= \xbt,$ $E: \xbt $

generated by $ \xCf ADE$ from $ \xbr_{2},$ $ \xbt $

prerequisite -

common part -

$ \xbr_{3}$ cannot be interpreted as a derived tripel, as there was a
prerequisite
used in its derivation (B), but the common part in $ \xbr_{3}$ is $ \xCQ
.$
 \xDH
$ \xbr_{4}:$ $A: \xbr_{3}= \xbs,$ $E: \xbr_{3}= \xbt,$ $B: \xbs = \xbt $

generated by $ \xCf AEB$ from $ \xbr_{3},$ $ \xbt $

prerequisites -

common part $B$

$ \xbr_{4}$ can be interpreted as the common part $B$ contains all
prerequisites
used in its derivation. $ \xCf ABE$ is the only non-trivial derived
tripel.

Note that we could, e.g., also have replaced ACD by $AC' (DC''),$ where
$C=C' \xcv C'',$
using rule (c), Weak Union.
 \xEJ

\ee

\be

$\hspace{0.01em}$


\label{Example R-4}

$BA(CD),$ $DF(CE),$ $(AB)(CD)(EF)$ $ \xch $ $B(ADF)(CE)$:

(See Diagram 
\ref{Diagram Rule-Const-1} (page 
\pageref{Diagram Rule-Const-1})  lower part.)

This example shows that we may need an assumption in the interior of the
tree (in the construction of $ \xbr_{3},$ we use $D: \xbs = \xbt).$
 \xEI
 \xDH
$ \xbr_{1}:$ $A: \xbs = \xbt,$ $B: \xbs,$ $C: \xbt,$ $D: \xbt $

generated by $BA(CD)$ from $ \xbs,$ $ \xbt $

prerequisites $ \xCf A$

common part $ \xCf A$
 \xDH
$ \xbr_{2}:$ $C: \xbt,$ $D: \xbs,$ $E: \xbt,$ $F: \xbs = \xbt $

generated by $DF(CE)$ from $ \xbs,$ $ \xbt $

prerequisite $F$

common part $F$
 \xDH
$ \xbr_{3}:$ A: $ \xbr_{1}= \xbs = \xbt,$ $B: \xbr_{3}= \xbs,$ $C:
\xbr_{1}= \xbr_{2}= \xbt,$ $D: \xbr_{1}= \xbr_{2}= \xbs = \xbt,$ $E:
\xbr_{2}= \xbt,$ $F: \xbr_{2}= \xbs = \xbt $

generated by $(AB)(CD)(EF)$ from $ \xbr_{1},$ $ \xbr_{2}$

prerequisite $D$

common part $ \xCf ADF$

So $ \xbr_{3}$ can be seen as the derived tripel $B(ADF)(CE)$ (but NOT as
$(AB)(DF)(CE)$
etc., as $ \xCf DF$ does not contain $ \xCf ADF.$
 \xEJ

\ee

\be

$\hspace{0.01em}$


\label{Example R-5}

$(AA')BC,$ $AD(CD'),$ $(AB')C(C' D),$ $(A' B')C(C' D'),$ $(AD)(B' CC'
)(A' D'),$
$BC(ADD')$ $ \xch $ $A(BD)(CD')$:

(See Diagram \ref{Diagram Rule-Const-3} (page \pageref{Diagram Rule-Const-3}).)

This example shows that we may need an equality
(here $ \xba $ and $ \xbb $ in the construction of $ \xbr_{5})$ which is
not
related to $ \xbs $ and $ \xbt.$ Of course, we cannot use it as an
assumption, but we
know the equality by construction.

$ \xba $ and $ \xbb $ will not be known, they are fixed, unknown
fragments.
 \xEI
 \xDH
$ \xbr_{1}:$ $A: \xbs,$ $A': \xbs,$ $B: \xbs = \xbt,$ $B': \xba,$
$C: \xbt $

generated by $(AA')BC$ from $ \xbs,$ $ \xbt $

prerequisites $B$

common part $B$
 \xDH
$ \xbr_{2}:$ $A: \xbs,$ $C: \xbt,$ $C': \xbb,$ $D: \xbs = \xbt,$ $D'
: \xbt $

generated by $AD(CD')$ from $ \xbs $ and $ \xbt $

prerequisite $D$

common part $D$
 \xDH
$ \xbr_{3}:$ $A: \xbs,$ $B': \xba,$ $C: \xbt,$ $C': \xbb,$ $D: \xbs
= \xbt $

generated by $(AB')C(C' D)$ from $ \xbr_{1}$ and $ \xbr_{2}$

prerequisite -

common part $D$
 \xDH
$ \xbr_{4}:$ $A': \xbs,$ $B': \xba,$ $C: \xbt,$ $C': \xbb,$ $D':
\xbt $

Generated by $(A' B')C(C' D')$ from $ \xbr_{1}$ and $ \xbr_{2}$

prerequisites -

common part -
 \xDH
$ \xbr_{5}:$ $A: \xbs,$ $A': \xbs,$ $B': \xba,$ $C: \xbt,$ $C':
\xbb,$ $D: \xbt,$ $D': \xbt $

generated by $(AD)(B' CC')(A' D')$ from $ \xbr_{3}$ and $ \xbr_{4}$

prerequisites - (note that equality on $B' $ and $C' $ is by construction
of $ \xbr_{3}$ and $ \xbr_{4},$ and not by a prerequisite on $ \xbs $ and
$ \xbt)$

common part: $D$
 \xDH
$ \xbr_{6}:$ $A: \xbs,$ $B: \xbs = \xbt,$ $C: \xbt,$ $D: \xbs = \xbt,$
$D': \xbt $

generated by $BC(ADD')$ from $ \xbr_{1}$ and $ \xbr_{5}$

prerequisites -

common part: $ \xCf BD$

Thus, $ \xbr_{6}$ may be seen as derived tripel $A(BD)(CD')$
 \xEJ

\ee

$ \xCO $


\begin{thebibliography}{xxxxxx}

\addcontentsline{toc}{section}{References}


\bibitem[Daw79]{Daw79}
A. P. Dawid, ``Conditional independence in statistical theory'', Journal
of the Royal Statistical Society, Series $B,$ $41(1):1-31,$ 1979

\bibitem[GS10]{GS10}
D. Gabbay, K. Schlechta,
``Conditionals and modularity in general logics'',
To appear (Springer, approx. spring 2011),
Preliminary version in arxiv.org

\bibitem[Par96]{Par96}
R. Parikh, ``Belief, belief revision, and splitting languages'',
Moss, Ginzburg and de Rijke (eds.)
Proceed. Logic, Language and Computation, CSLI 1999, pp. 266--278

\bibitem[Pea88]{Pea88}
J. Pearl, ``Probabilistic Reasoning in Intelligent Systems'',
Morgan Kaufmann, San Mateo, Cal., 1988

\bibitem[Spo80]{Spo80}
W. Spohn, ``Stochastic independence, causal independence, and
shieldability'', Journal of Philosophical Logic 9 (1980) 73-99

\end{thebibliography}
\end{document}